\documentclass[10pt,a4paper]{article}

\usepackage{amsmath}
\usepackage{graphicx,latexsym,euscript,makeidx,color}
\usepackage{amsmath,amsfonts,amssymb,amsthm,mathrsfs}
\usepackage{enumerate}

\usepackage{amssymb}
\usepackage{amsmath}
\usepackage{amsfonts}
\usepackage{graphicx}
\usepackage{graphicx}          
\usepackage{color}

\usepackage[colorlinks,linkcolor=blue,anchorcolor=green,citecolor=red]{hyperref}

\usepackage{geometry}
\font\tenbb=msbm10 \font\sevenbb=msbm7 \font\fivebb=msbm5
\newfam\bbfam
\scriptscriptfont\bbfam=\fivebb \textfont\bbfam=\tenbb
\scriptfont\bbfam=\sevenbb

\usepackage{graphics} 
\usepackage{epsfig} 
\usepackage{mathptmx} 
\usepackage{times} 
\usepackage{amsmath} 
\usepackage{amssymb}  
\usepackage{color}
\usepackage{graphicx,colordvi}

\newtheorem{theorem}{\indent Theorem}[section]
\newtheorem{definition}[theorem]{\indent Definition}
\newtheorem{proposition}[theorem]{\indent Proposition}

\newtheorem{lemma}[theorem]{\indent Lemma}
\newtheorem{remark}[theorem]{\indent Remark}

\makeatletter
   \renewcommand{\theequation}{%
            \thesection.\arabic{equation}}
   \@addtoreset{equation}{section}
\makeatother

\allowdisplaybreaks[4]

\begin{document}

\title{\bf Mean Field Games for Multi-agent Systems with Multiplicative Noises
\thanks{This work was supported in part by the National Natural Science Foundation of China under Grants 61773241, 61773222, 61573221 and 61633014, and in part by Taishan Scholar Construction Engineering by Shandong Government.}
}
\author{Bing-Chang Wang\thanks{School of Control Science and Engineering,
					Shandong University, Jinan 250061, P. R. China. Email: {\tt bcwang@sdu.edu.cn}.}~~~~~~
Yuan-Hua Ni\thanks{College of Artificial Intelligence, Nankai University, Tianjin 300350, P. R. China. Email: {\tt yhni@nankai.edu.cn}.}~~~~~~
Huanshui Zhang\thanks{School of Control Science and Engineering, Shandong University, Jinan 250061, P. R. China. Email: {\tt hszhang@sdu.edu.cn}.}}
\maketitle

{\bf Abstract:} This paper studies mean field games for multi-agent systems with control-dependent multiplicative noises. For the general systems with nonuniform agents, we obtain a set of decentralized strategies by solving an auxiliary limiting optimal control problem subject to consistent mean field approximations.
 The set of decentralized strategies is further shown to be an $\varepsilon$-Nash equilibrium. For the integrator multiagent systems, we design a set of $\varepsilon$-Nash strategies by exploiting the convexity property of the limiting problem. It is shown that under the mild conditions all the agents achieve mean-square consensus.



\section{Introduction}
The methodology of mean field games was proposed to investigate decentralized strategies in large population dynamic games  \cite{HCM03}, \cite{HCM07}, \cite{HMC06}, \cite{LL07}; it is effective to design asymptotic Nash equilibria for many players which are individually negligible but collectively affect a particular player. By identifying a consistency relationship between the individual's optimal response  and the population macroscopic behavior, one may obtain a  fixed-point equation to specify the mean field effect.  This procedure leads to a set of  decentralized  strategies as an $\epsilon$-Nash equilibrium for the actual model with large population.
For some aspects of mean field games, the readers are referred to the linear-quadratic (LQ, for short) framework \cite{HCM03}, \cite{HCM07}, \cite{TZB14}, \cite{NLZ15}, \cite{WZ12}, mean field games with a major player \cite{BFY13}, \cite{huang2010large}, \cite{wang2012major}, \cite{WZ14}, the nonlinear mean field games \cite{CD13}, \cite{HMC06}, oblivious equilibria for large-scale Markov decision processes \cite{weintraub2008markov}, and mean field games with Markov jump parameters \cite{WZ12}, \cite{WZ13}. For the survey on mean field games, see, e.g., \cite{BFY13}, \cite{GS13}, and \cite{CHM17}.

Concerned with mean field games for stochastic systems, most existing literature focuses on the case with additive noise, namely, the intensity of noise is independent of the system state. Sometimes, such kind of noises is not sufficient to model the practical situations.
%
%
Alternatively, multiplicative noise is another realistic description for uncertainty; in this case, the noise intensity depends on the system state. Stochastic systems with multiplicative noises have been extensively studied in the past half century in the control and mathematical communities. From the viewpoint of mathematics, almost all the theories about stochastic differential equations (SDEs, for short) are for the case with multiplicative noises, and there are lots of practical motivations to study such kinds of SDEs \cite{O11}.
The study in control community is also popular; to mention a few, see monographs \cite{FR75}, \cite{YZ99} and some recent papers on stochastic control \cite{ZX17}, \cite{ZZC08}.

In particular, some researchers have paid attention to multi-agent systems with multiplicative noises, that is, the intensity of noise depends on agents'
states.
Multiplicative noises in multi-agent systems can be generated by data transmission channels, during the
propagation of radio signals or under signal processing by receivers or detectors \cite{WE13}.
In \cite{NL13}, the authors investigated the consensus problems for the continuous-time
systems with multiplicative noises and the noise intensities are proportional to the
absolute value of the relative states of agents and their neighbors.
Then this work was extended to the discrete-time version in \cite{Long2015}. Further, the work \cite{LZ14} studied the distributed consensus with the general multiplicative noises and developed some small consensus gain theorems to give sufficient conditions for mean square and almost sure consensus under undirected topologies. Moreover, the authors in \cite{Zong-1} investigated stochastic consentability of linear multi-agent systems with time delays and multiplicative noises; and \cite{Zong-2} developed the consensus conditions of continuous-time multi-agent systems with additive and multiplicative measurement noises.
Note that the methodology of handling the multiplicative noises of multi-agent systems is different from the stochastic-approximation-type algorithms \cite{HM09}, \cite{LZ09}, \cite{WZ09}, which are powerful tools to deal with multi-agent systems with additive noisy measurements.
%
%

Due to its wide popularity and practical applications \cite{CD18, GL11, HN16}, multiplicative noise should be taken care of when considering mean field games, which is the theme of this paper.
An LQ mean field game with multiplicative noises is introduced in Section \ref{sec2.3.1}, where the diffusion terms include both the state and the control variables.
Due to the appearance of multiplicative noise, the strategy design and analysis of this LQ mean field game is not a trivial extension of those with additive noises \cite{ZX17}, \cite{LFZ17}.
For instance, a stochastic linear system with additive noise cannot attain asymptotic stability (in the normal sense), while a stochastic system with multiplicative noise could be asymptotically stable.
Furthermore, the optimal control law for a stochastic system with additive noise and a quadratic cost is no longer optimal for the case with multiplicative noise.
Indeed, research on stochastic systems
with additive noise shares similar methods and results to deterministic systems. In
contrast, it is more complicate to study stochastic control with multiplicative noise,
particularly with control-dependent noise (e.g., stochastic maximum principle) \cite{YZ99}.
In Section \ref{Sec:Strategy-design}, we start by solving an auxiliary limiting optimal control problem and next construct an equation system of mean field effect by the consistent mean field approximation. For the case of the uniform agents, the conditions for existence and uniqueness of the solution to the consistency equations
  are obtained. With the help of the mean field effect and the solution to the limiting problem, we design a set of decentralized strategies, which is further shown to be an $\varepsilon$-Nash equilibrium.
%

In Section \ref{Sec:multi-agent},  we investigate single- and multiple- integrator multi-agent systems with multiplicative noises. The dynamics of each agent is a integrator system with control-dependent multiplicative noise, and the control weight in the cost function is not limited to be positive. By exploiting the coercivity property of the limiting problem, the proposed strategies is proved to be an $\varepsilon$-Nash equilibrium. Particularly, we give the explicit conditions to ensure all the agents achieve mean-square consensus. The closely related works \cite{NCMH13}, \cite{NCM14} investigated the consensus problem for single-integrator systems with additive noises from the optimization principle, and under the proposed strategies the mean of agents' states asymptotically achieve consensus, whereas the own states of agents can not attach consensus.
In Section \ref{Sec:example}, we give a numerical example to illustrate the effectiveness of the proposed consensus strategies for single-integrator multi-agent systems with multiplicative noises. Furthermore, we compare the consensus effect between this model and multi-agent systems with additive noises. 
It is shown that agents with multiplicative noises
can reach consensus, while agents with additive noises do not achieve consensus exactly.
Section \ref{Sec:conclude} concludes the paper.

Compared with previous works, the main contributions of the paper are summarized as follows:
(i) We study mean field games with control-dependent multiplicative noises and obtain decentralized $\varepsilon$-Nash strategies. 
Different from the case with additive noise, we adopt the Lyapunov function method to show the uniform stability of closed-loop system, 
and then obtain the asymptotic equilibrium property of decentralized strategies.
(ii) For the case with uniform agents, we give an easy-verified condition to ensure the existence and uniqueness of a solution to consistency equations. (iii) We further investigate (single- and multiple-) integrator multi-agent systems with multiplicative noises, where the control weight in the cost is not limited to be positive. By exploiting the convexity and coercivity properties of the limiting problem, it is shown that the set of decentralized strategies is an $\varepsilon$-Nash equilibrium and all the agents of closed-loop system converge to a constant, which is the limit of the initial state average of all agents as $N\to\infty$.


The following notation will be used throughout this paper. $\|\cdot\|$
denotes the Euclidean vector norm or Peron-Frobenious matrix norm. For a vector $z$ and a matrix $Q$, $\|z\|_Q^2= z^TQz$, and $Q>0$ (or $Q\geq0$) means that $Q$ is positive definite (or semi-positive definite).
$C([0,\infty),\mathbb{R}^n)$ is the space of all the $n$-dimensional continuous functions on $[0,\infty)$, and $C_{\rho/2}([0,\infty), \mathbb{R}^n)$ is a subspace of $C([0,\infty),\mathbb{R}^n)$ which is given by $\{f|\int_0^{\infty}e^{-\rho t}\|f(t)\|^2dt<\infty  \}.$
 For a family of $\mathbb{R}^n$-valued random variables $\{x(\lambda),
 \lambda\geq 0\}$,
$\sigma(x(\lambda), \lambda\leq t)$ is the $\sigma$-algebra  generated by the collection of random variables.
For two sequences $\{a_n, n = 0, 1, \cdots\}$ and $\{b_n, n = 0, 1,  \cdots\}$, $a_n = O(b_n)$ denotes
$\limsup_{n\to\infty}|{a_n}/{b_n}|\leq c$, and $a_n = o(b_n)$ denotes
$\limsup_{n\to\infty}|{a_n}/{b_n}|=0$.
 For convenience of  presentation, we use $c, c_1,c_2,\cdots$ to
denote generic positive constants, which may vary from place to place.

\section{Problem Formulation}\label{sec2.3.1}

Consider a large population system with $N$ agents. For $i=1,...,N$, agent $i$ evolves according to the following SDE:
\begin{equation}\label{eq1}
\begin{aligned}
dx_i(t) = [A(\theta_i)x_i(t)+Bu_i(t)]dt+[Cx_i(t)+Du_i(t)]dW_i(t),
\end{aligned}
\end{equation}
where $x_i\in \mathbb{R}^n$ and $u_i\in\mathbb{R}^r$ are the state and input.  The underlying filtered probability space is $(\Omega, \{{\mathcal F}_t\}_{t\geq 0}, \mathbb{P})$, where $\{{\mathcal F}_t\}_{t\geq 0}$ is a collection of non-decreasing $\sigma$-algebras.
$W_i=\{W_i(t), t\geq 0\}$ is a $1$-dimensional standard Brownian motion, which is adapted to $\{\mathcal{F}_t\}_{t\geq 0}$, and $W_i, i= 1,...,N,$ are assumed to be mutually independent.
The matrices $A(\cdot), B, C$ and $D$ are deterministic and have compatible dimensions. Here, the parameters $\theta_i, 1\leq i\leq N$, are to model this population of nonuniform agents. For notational simplicity, the time argument of a process ($x_i,u_i$, etc.) is sometimes suppressed, and $A(\theta_i)$ is denoted as $A_i$ or $A_{\theta_i}$. 

The cost function of agent $i$ is given by
\begin{equation}\label{eq2}
\begin{aligned}
J_i(u_i,u_{-i})=\mathbb{E}\int_0^{\infty}e^{-\rho t} \Big\{ \big\| x_i(t) -
	x^{(N)}(t) \big \|^2_{Q}+\|u_i(t)\|^2_{R} \Big \}dt,
\end{aligned}
\end{equation}
where $\rho\geq 0$, and $Q, R$ are symmetric matrices with compatible dimensions.
%
%
%
$x^{(N)}(t)=\frac{1}{N}\sum_{j=1}^Nx_j(t)$, and $u_{-i}=\{u_1, \ldots,u_{i-1},u_{i+1}, \ldots, u_N\}$.
The admissible control set of agent $i$ is given by
$$
\begin{aligned}
{\mathcal U}_{i} &=\Big\{u_i\, \big{|}\,
u_i(t)\in {\mathcal F}_t^{i},\ \mathbb{E}\int_0^{\infty}e^{-\rho t}\|x_i(t)\|^2dt<\infty\Big\}, 
\end{aligned}
$$
where ${\mathcal F}_t^{i}=\sigma(x_i(s),s\leq t)$, $i=1,\cdots,N$.
Denote $\mathcal{U}_d\triangleq \mathcal{U}_1 \times \cdots \times \mathcal{U}_{N}$.
For comparison to $\mathcal{U}_d$, we introduce the centralized admissible control set:
$$\begin{aligned}
{\mathcal U}_{c} &=\Big\{u_i|\
u_i(t)\in \sigma(\bigcup_{i=1}^N{\mathcal F}_t^{i}), \mathbb{E}\int_0^{\infty}e^{-\rho t}\|x_i(t)\|^2dt<\infty \Big\}.
\end{aligned}
$$

In this paper, we will study the following problem.

\textbf{Problem (P).} \emph{Seek a set of control strategies within
	$\mathcal{U}_d$
	that minimize $J_i(u_i,u_{-i}), i=1,...,N$ in a game-theoretic setting. Specifically, find an $\varepsilon$-Nash equilibrium in the decentralized control set $\mathcal{U}_d$, where $\varepsilon\to 0$ as $N \to \infty$.}

For a large population, a natural way of modeling parameters $\theta_1,\cdots,\theta_N$ is to view them as a sequence of observations on the same underlying distribution function
such that this sequence exhibits certain statistical properties when $N\to\infty$. For this, assume that $\theta_1,\cdots,\theta_N$ are valued in a compact set $\Theta \in\mathbb{R}^d$. Define the associated empirical distribution function $F_N(\theta)=\frac{1}{N}\sum_{i=1}^NI_{[\theta_i\leq \theta]}$, where $I_{[\theta_i\leq \theta]}=1$ if $\theta_i\leq \theta$ and $I_{[\theta_i\leq \theta]}=0$ otherwise.

We make the following assumptions:

\textbf{A1)} $x_i(0), i=1,...,N$, are mutually independent and have the same mathematical expectation. $\mathbb{E}x_i(0)\equiv\bar{x}_0$, and there exists a constant $c_0$ (independent of $N$) such that $max_{1\leq i \leq N}\mathbb{E}\|x_i(0)\|^2<c_0$. Furthermore, $\{x_i(0), i=1,...,N\}$ and
$\{W_i, i=1,...,N\}$ are independent of each other.

\textbf{A2)} $A(\theta)$ is a continuous matrix function of $\theta\in \Theta$, where $\Theta$ is a compact subset of $\mathbb{R}^d$. There exists a distribution function $F(\cdot)$ such that $F_N$ converges weakly to $F$, i.e., for any bounded continuous function $\psi(\theta)$ on $\Theta$, $lim_{N\to\infty}\int_\Theta\psi(\theta) dF_{N}(\theta)=\int_\Theta\psi(\theta) dF(\theta)$.

\textbf{A3)} For any $i\in \{1,...,N\}$, the system (\ref{eq1}) is $\rho$-stabilizable\footnote{A similar definition is $\rho$-stability. A system $$dx(t)=Ax(t)dt+Cx(t)dW(t)$$ is said to be $\rho$-stable (or [A,C] is $\rho$-stable, for short)
	if $\mathbb{E}\int_0^{\infty}e^{-\rho t}\|x(t)\|^2dt<\infty$. See \cite{ZZC08} for further results about $\rho$-stability and $\rho$-stabilizability.}
in the sense that there exists a feedback control law $u_i=K_ix_i$ such that
for any initial state $x_i(0)\in \mathbb{R}^n$, the state of closed-loop system
$$dx_i(t)=(A_i+BK_i)x_i(t)dt+(C+DK_i)x_i(t)dW_i(t)$$
satisfies
$$\mathbb{E}\int_0^{\infty}e^{-\rho t}\|x_i(t)\|^2dt<\infty.$$


\textbf{A4)} $R>0$, $Q\geq0$. The system
$$\left\{\begin{aligned}
dz_i(t)=&(A_i-\frac{\rho}{2}I)z_i(t)dt+Cz_i(t)dW_i(t),\\
y_i(t)=&\sqrt{Q}z_i(t)
\end{aligned}
\right.
$$
(or simply $[A_i-\frac{\rho}{2}I, C,\sqrt{Q}]$) is exactly detectable, namely, if $y_i(t)=0, \hbox{a.s.},\ t\in [0,T], \forall T>0$ implies
$lim_{t\to\infty} \mathbb{E}\|z_i(t)\|^2\to 0$.

\begin{remark}\label{rem1}
\Blue{Assumption \textbf{A3)} is necessary. Otherwise, the cost function cannot be finite, to say nothing of the optimal control. \textbf{A4)} could be relaxed. Actually, provided some convexity condition is satisfied, neither $R$ nor $Q$ is limited to be positive definite. For instance, $R$ is allowed to be negative in Section 4. In \cite{WH2015},
$Q=\left[
\begin{array}{cc}
  0 &-0.5 \\
  -0.5& 0
\end{array}\right]$ is not semi-positive definite.}
\end{remark}

\begin{remark}
 By \cite[Theorem 1]{RZ00} The system (\ref{eq1}) is $\rho$-stabilizable
if and only if there exists a matrix $K$ and a positive definite matrix $X_i$ such that
$$\rho X_i =(A_i+BK)^TX_i+X_i(A_i+BK)+(C+DK)^TX_i(C+DK)<0.$$
By \cite[Theorem 3.1]{ZZC08}, the system (\ref{eq1}) is exactly detectable
if and only if  there does not exist nonzero symmetric matrix $Z$ such that
$$AZ+ZA^T+CZC^T=\lambda Z, \ QZ=0, \ Re\lambda\geq 0. $$
See more equivalent conditions in verifying $\rho$-stabilizability and exact detectability in \cite{RZ00, ZZC08}.
\end{remark}

\section{Strategy Design and Analysis}\label{Sec:Strategy-design}

\subsection{Optimal control for the limiting problem}

%
Assume that $\bar{x}(\cdot) \in C_{\rho/2}([0,\infty), \mathbb{R}^n)$ is given for approximation of $x^{(N)}(\cdot)$. Replacing $x^{(N)}(\cdot)$ in (\ref{eq2}) by $\bar{x}(\cdot)$, we have a standard optimal tracking problem.

\textbf{Problem (P1)}. \emph{For $i=1,...,N$, find $\bar{u}_i\in \mathcal{U}_i$ such that
	\begin{equation*}
		\bar{J}_i(\bar{u}_i)=\min_{u_i\in \mathcal{U}_i}\bar{J}_i(u_i)
	\end{equation*}
	holds, where $\bar{J}_i(u_i)$ is given by
	\begin{equation*}
		\bar{J}_i(u_i)=\mathbb{E}\int_0^{\infty}e^{-\rho t}\Big\{\big\|x_i(t)
		- \bar{x}(t)\big\|^2_{Q}
		+\|u_i(t)\|^2_{R}\Big\}dt.
	\end{equation*}
}

For a general initial state $x_i(t)=x_i\in \mathbb{R}^n$, define the value function
%
%
\begin{equation*}
\begin{aligned}
V_i(t,x_i) = \inf_{u_i\in{\mathcal U}_{i}} \mathbb{E} \left[  \int_{t}^{\infty}e^{-\rho(\tau-t)}\Big\{\big\|x_i(\tau)-\bar{x}(\tau) \big\|^2_{Q}+\big\|u_i(\tau)\big\|^2_{R}\Big\}d\tau\Big|x_i(t)=x_i  \right].
\end{aligned}
\end{equation*}
We introduce the HJB equation:
\begin{eqnarray}\label{eq6n}
&&\hspace{-3em}\rho V_i(t,x_i) = \inf_{u_i\in \mathbb{R}^r}
	\Big{\{} \frac{\partial V_i(t,x_i)}{\partial t} +\frac{\partial V_i^T(t,x_i)}{\partial x_i}\left(A_ix_i+Bu_i\right)
	 +\frac{1}{2}(Cx_i+Du_i)^T\frac{\partial^2 V_i(t,x_i)}{\partial x^2 _i}(Cx_i+Du_i)\nonumber\\
&&\hspace{-3em}\hphantom{\rho V_i(t,x_i) =}	 +\|x_i - \bar{x}(t)\|^2_{Q} + \| u_i\| ^2_{R} \Big{\}}.
\end{eqnarray}
Letting $V_i(t,x_i)=x_i^TP_ix_i+2s_i^T(t)x_i+g_i(t)$, then the optimal control is given by
\begin{align}\label{eq7n}
\bar{u}_i=&-\frac{1}{2}(R+D^TP_iD)^{-1}\big[B^T\frac{\partial V_i}{\partial x_i}+D^T\frac{\partial ^2V_i}{\partial x_i^2}Cx_i\big]\cr
=&-(R+D^TP_iD)^{-1}[(B^TP_i+D^TP_iC) x_i+B^Ts_i].
\end{align}
Substituting (\ref{eq7n}) into (\ref{eq6n}), we have
\begin{eqnarray*}
&&\hspace{-3em}\rho(x_i^TP_ix_i+2s_i^Tx_i+g_i)=x_i^T[A^T_iP_i+P_iA_i+C^TP_iC+Q-(B^TP_i+D^TP_iC)^T(R+D^TP_iD)^{-1}(B^TP_i+D^TP_iC)]x_i\cr
&&\hspace{-3em}\hphantom{\rho(x_i^TP_ix_i+2s_i^Tx_i+g_i)=}+2\left[\frac{ds_i}{dt}+[A_i-B(R+D^TP_iD)^{-1}(B^TP_i+D^	TP_iC)]^Ts_i-Q\bar{x}\right]x_i\cr
&&\hspace{-3em}\hphantom{\rho(x_i^TP_ix_i+2s_i^Tx_i+g_i)=}+\frac{dg_i}{dt}-s_i^TB(R+D^TP_iD)^{-1}B^Ts_i+\bar{x}^TQ\bar{x}.
\end{eqnarray*}
This yields
\begin{align}\label{eq6}
\rho P_i =&A^T_iP_i+P_iA_i+C^TP_iC+Q-(B^TP_i+D^TP_iC)^T(R+D^TP_iD)^{-1}(B^TP_i+D^TP_iC),\\\label{eq7a}
\rho s_i=&\frac{ds_i}{dt}+[A_i-B(R+D^TP_iD)^{-1}(B^TP_i+D^TP_iC)]^Ts_i-Q\bar{x},\\
\rho g_i=&\frac{dg_i}{dt}-s_i^TB(R+D^TP_iD)^{-1}B^Ts_i+\bar{x}^TQ\bar{x}. \label{eq8}
\end{align}
Here, $s_i(0)$ and $g_i(0)$ are not pre-specified and to be determined by  $s_i\in {C}_{\rho/2}([0,\infty), \mathbb{R}^n)$ and $g_i\in C_{\rho}([0,\infty), \mathbb{R}^n)$, respectively.

\Blue{We now provide a preliminary result for LQ control with multiplicative noise.}
\begin{lemma} \label{thm1}
	Assume that \textbf{A1)}-\textbf{A4)} hold and $\bar{x}\in C_{\rho/2}([0,\infty), \mathbb{R}^n)$ is given. For Problem (P1), the following statements hold.
	
	i) The algebraic Riccati equation (\ref{eq6}) admits a unique semi-positive definite solution, and the closed-loop system
	\begin{equation}\label{eq10a}
	dx_i(t)=\bar{A}_ix_i(t)dt+\bar{C}x_i(t)dW_i(t)
	\end{equation}
	or $[\bar{A}_i,\bar{C}]$ is $\rho$-stable, where $\bar{A}_i\stackrel{\Delta}{=}A_i-B(R+D^TP_iD)^{-1}(B^TP_i+D^TP_iC)$ and $\bar{C}\stackrel{\Delta}{=}C-D(R+D^TP_iD)^{-1}(B^TP_i+D^TP_iC)$, $1\leq i\leq N$.

	ii) The equation (\ref{eq7a}) admits a unique solution $s_i\in {C}_{\rho/2}([0,\infty), \mathbb{R}^n)$, and (\ref{eq8}) admits a unique solution $g_i\in C_{\rho}([0,\infty), \mathbb{R})$.

	iii) Problem (P1) admits a unique optimal control
	\begin{equation}\label{optimal-L}
	\bar{u}_i(t)=-(R+D^TP_iD)^{-1}[(B^TP_i+D^TP_iC) x_i(t)+B^Ts_i(t)], \ t\geq 0.
	\end{equation}
	and the optimal cost is given by
	\begin{equation*}\label{eq10n}
		\inf_{u_i\in {\mathcal U}_{i}} \bar{J}_i(u_i)=
		\mathbb{E}[x^T_i(0)P_ix_i(0)]+2s_i^T(0)\bar{x}_0+g_i(0).
	\end{equation*}
	
\end{lemma}

\emph{Proof.}
See Appendix \ref{app1}.
\hfill$\Box$

\subsection{Design of control strategies}

Following the standard approach in mean field games \cite{HCM07}, \cite{LZ08}, we construct the consistency equations
\begin{align}\label{eq9a}
\rho s_{\theta}(t)=& \frac{ds_{\theta}}{dt}+
\bar{A}_{\theta}^Ts_{\theta}(t)-Q\bar{x}(t),\\ \label{eq10}
\frac{d\bar{x}_{\theta}}{dt}=&\bar{A}_{\theta}\bar{x}_{\theta}(t)-B(R+D^TP_{\theta}D)^{-1}B^Ts_{\theta}(t),\\ \label{eq11a}
\bar{x}(t)=&\int_{\Theta}\bar{x}_{\theta}(t)dF(\theta)
\end{align}
where $\bar{A}_{\theta}=A_{\theta}-B(R+D^TP_{\theta}D)^{-1}(B^TP_{\theta}+D^TP_{\theta}C)$, $\bar{x}_{\theta}(0)=\bar{x}_0$, and $s_{\theta}(0)$ is to be determined.
In the above, $\bar{x}_{\theta}$ is regarded as the expectation of the state given the parameter $\theta$ in the individual dynamics. {The last equation is due to the consistency requirement for the mean field approximation.}

\Blue{The existence and uniqueness of a solution to (\ref{eq9a})-(\ref{eq11a}) may be obtained by using fixed-point methods similar to those in \cite{HCM07} and \cite{WZ12}.}
For further analysis,
we make the following assumption.

\textbf{A5)} There exists a solution
$\bar{x}\in C_{\rho/2}([0,\infty), \mathbb{R}^n)$ to (\ref{eq9a})-(\ref{eq11a}).

For a population-size system in Problem (P), we may design the control strategy as follows: 
\begin{equation}\label{eq10bb}
\begin{aligned}
\hat{u}_i(t)=-(R+D^TP_iD)^{-1}[(B^TP_i+D^TP_iC) x_i(t)+B^Ts_i(t)], \ t\geq 0,\ i=1,\cdots, N,
\end{aligned}
\end{equation}
where $s_i$ is the solution of equation system (\ref{eq9a})-(\ref{eq11a}).

\subsection{The case with uniform agents}

In what follows, we consider the consistency equations for the case with uniform agents, namely, $A({\theta_i})\equiv A$ in (\ref{eq1}) is independent of $\theta_i, i=1,\cdots,N$. In this case, the equations (\ref{eq9a})-(\ref{eq11a}) read as
\begin{align}\label{eq19a}
\rho s(t)=& \frac{ds}{dt}+
\bar{A}^Ts(t)-Q\bar{x}(t),\\ \label{eq19b}
\frac{d\bar{x}}{dt}=&\bar{A}\bar{x}(t)-B(R+D^TPD)^{-1}B^Ts(t),
\end{align}
where $\bar{A}{=}A-B(R+D^TPD)^{-1}(B^TP+D^TPC).$ Denote
$$M=\left[
\begin{array}{cc}
\bar{A}-\frac{\rho}{2} I& B(R+D^TPD)^{-1}B^T\\
-Q&-\bar{A}^T+\frac{\rho}{2} I
\end{array}
\right].$$
\Blue{
We now give an easy-verified condition to guarantee \textbf{A5)}.
\begin{theorem}\label{thm3.2}
	Let \textbf{A3)} and \textbf{A4)} hold. If 
the real part of any eigenvalue of $M$ is not zero (no eigenvalue of $M$ is on the imaginary axis), then (\ref{eq19a})-(\ref{eq19b}) have a unique solution $(s,\bar{x})\in C_{\rho/2}([0,\infty), \mathbb{R}^{2n})$.
\end{theorem}}

\emph{Proof.} Letting $s(t)=K\bar{x}(t)+\varphi(t)$, $t\geq0$, from (\ref{eq19b}) and (\ref{eq19a}) we have
\begin{align*}
\frac{ds}{dt}=&K \frac{d\bar{x}}{dt}+ \frac{d\varphi}{dt}\cr
=&K[\bar{A}\bar{x}-B(R+D^TPD)^{-1}B^T(K\bar{x}+\varphi)]+ \frac{d\varphi}{dt} \cr
=&(\rho I-
\bar{A}^T)(K\bar{x}+\varphi)+Q\bar{x},
\end{align*}
which implies
\begin{align}\label{eq24}
K(\bar{A}-\frac{\rho}{2}I)&+(\bar{A}-\frac{\rho}{2}I)^TK -KB(R+D^TPD)^{-1}B^TK-Q=0,\\
\label{eq25}
\frac{d\varphi}{dt}+\big[\bar{A}-&\rho I-B(R+D^TPD)^{-1}B^TK\big]^T \varphi=0.
\end{align}
From Lemma \ref{thm1}, $\bar{A}_i-\frac{\rho}{2}I$
is Hurwitz. 
Since no eigenvalue of $M$ is on the imaginary axis, then by \cite{M77} we have that (\ref{eq24}) admits a unique solution such that $\bar{A}-\frac{\rho}{2} I-B(R+D^TPD)^{-1}B^TK$ is Hurwitz. This implies that each eigenvalue of $\bar{A}-{\rho} I-B(R+D^TPD)^{-1}B^TK$ has a real part lesser than $-\frac{\rho}{2}$. We may identify a unique $\varphi(0)=0$ such that $\varphi(t)\in C_{\rho/2}([0,\infty), \mathbb{R}^n)$ and the resulting solution is $\varphi=0$.
By (\ref{eq19b}), we have
\begin{equation}\label{eq26}
\frac{d\bar{x}}{dt}=\big[\bar{A}-B(R+D^TPD)^{-1}B^TK\big]\bar{x}(t), \quad \bar{x}(0)=\bar{x}_0.
\end{equation}
Since $\bar{A}-\frac{\rho}{2} I-B(R+D^TPD)^{-1}B^TK$ is Hurwitz, then all the eigenvalues of $\bar{A}-B(R+D^TPD)^{-1}B^TK$ have real parts less than $\frac{\rho}{2}$. Hence, $\bar{x}\in C_{\rho/2}([0,\infty), \mathbb{R}^n)$, and $s=K\bar{x}\in C_{\rho/2}([0,\infty), \mathbb{R}^n)$.

We now show the uniqueness. Suppose that $(s^{\prime},\bar{x}^{\prime})$ is another solution (in $ C_{\rho/2}([0,\infty), \mathbb{R}^{2n})$). Denoting $\varphi^{\prime}=s^{\prime}-K\bar{x}^{\prime}$, then $\varphi^{\prime}\in  C_{\rho/2}([0,\infty), \mathbb{R}^{n})$ satisfies ({\ref{eq25}}), which leads to $\varphi^{\prime}=0$. Subsequently, $\bar{x}^{\prime}$ satisfies ({\ref{eq26}}). Hence, $\bar{x}^{\prime}=\bar{x}$, which further implies $s^{\prime}=s$. \hfill$\Box$

\subsection{Asymptotic Nash equilibrium}
We first present a lemma on the $\rho$-stability of the following system
\begin{equation}\label{10a}
dx(t)=[Ax(t)+f(t)]dt+[Cx(t)+\sigma(t)]dW(t), \quad t\geq0,
\end{equation}
where $\{W(t), t\geq 0\}$ is a $1$-dimensional standard Brownian motion. A similar result was given in \cite{SY16}.
\begin{lemma}\label{lem1}
	Suppose that the system 
	$[{A}, {C}]$
	is $\rho$-stable, and $\mathbb{E}\|x(0)\|^2<\infty$. Then, {for any $f, \sigma\in C_{\rho/2}([0,\infty),\mathbb{R}^n)$ we have}
	\begin{equation*}
	\begin{aligned}
	 \mathbb{E}\int_0^{\infty}e^{-\rho t} \|x(t)\|^2dt
\leq  c\mathbb{E}\|x(0)\|^2 +c\mathbb{E}\int_0^{\infty}e^{-\rho t}[\|f(t)\|^2+\|\sigma(t)\|^2]dt< \infty.
	\end{aligned}
	\end{equation*}
\end{lemma}

\emph{Proof.}  Since  $[{A}, {C}]$ is $\rho$-stable, by the proof of i) of Lemma \ref{thm1} and \cite{HLY14}, there exists a matrix $P>0$ such that
$$P(A-\frac{\rho}{2}I)+(A-\frac{\rho}{2}I)^TP+C^TPC=-M<0.$$
Let $V(t)=e^{-\rho t}{x}^T(t)P{x}(t)$. By using It\^{o}'s formula, one has
\begin{eqnarray*}
&&\hspace{-2em}\mathbb{E}V(T)-\mathbb{E}V(0)
=\mathbb{E}\int_0^T e^{-\rho t}\big[{x}^T(t)(-\rho P+A^TP+PA+C^TPC){x}(t)\nonumber\\
&&\hspace{-2em}\hphantom{\mathbb{E}V(T)-\mathbb{E}V(0)
=}+2(Pf(t)+C^TP\sigma(t))^Tx(t)+\sigma^T(t)P\sigma(t)\big]dt\\
&&\hspace{-2em}\hphantom{\mathbb{E}V(T)-\mathbb{E}V(0)}=\mathbb{E}\int_0^Te^{-\rho t}\big[-x^T(t)Mx(t)+2(Pf(t)+C^TP\sigma(t))^Tx(t)+\sigma(t)^TP\sigma(t)\big]dt\\
&&\hphantom{\mathbb{E}V(T)-\mathbb{E}V(0)}\leq -c_1 \mathbb{E}\int_0^TV(t)dt+\mathbb{E}\int_0^T e^{-\rho t}[2(Pf(t)+C^TP\sigma(t))^Tx(t)+\lambda_{max}(P)\|\sigma(t)\|^2 ]dt\\
&&\hspace{-2em}\hphantom{\mathbb{E}V(T)-\mathbb{E}V(0)}\leq -\frac{c_1}{2} \mathbb{E}\int_0^TV(t)dt+c_2 \mathbb{E}\int_0^T e^{-\rho t}(\|f(t)\|^2+\|\sigma(t)\|^2)dt,
\end{eqnarray*}
where
$$c_1=\frac{\lambda_{min}(M)}{\lambda_{max}(P)},\quad  c_2>0.
$$
Here $\lambda_{min}(M)$ is the minimum eigenvalue of $M$ and $\lambda_{max}(P)$ is the maximum eigenvalue of $P$.
Then for any $T>0$, we obtain
\begin{equation*}
\begin{aligned}
\mathbb{E}\int_0^TV(t)dt &\leq\frac{2}{c_1}\Big[\mathbb{E}V(0)-\mathbb{E}V(T)+c_2\mathbb{E}\int_0^T e^{-\rho t}(\|f(t)\|^2+\|\sigma(t)\|^2)dt\Big]\\
&\leq\frac{2}{c_1}\Big[\lambda_{max}(P)\mathbb{E}\|x(0)\|^2+c_2\mathbb{E}\int_0^T e^{-\rho t}(\|f(t)\|^2 +\|\sigma(t)\|^2)dt\Big].\\
\end{aligned}	
\end{equation*}
Letting $T\to \infty$, this completes the proof.   \hfill$\Box$

After the strategy (\ref{eq10bb}) is applied, the closed-loop dynamics for agent $i$ may be written as follows:
\begin{equation*}
\begin{aligned}
d\hat{x}_i(t)= \ & \bar{A}_i\hat{x}_i(t)dt-B(R+D^TP_iD)^{-1}B^Ts_i(t)dt\\
&+[\bar{C}\hat{x}_i(t)-D(R+D^TP_iD)^{-1}B^Ts_i(t)]dW_i(t),\  1\leq i\leq N.\\
\end{aligned}
\end{equation*}

\begin{theorem}\label{thm2}
	Let \textbf{A1)}-\textbf{A5)} hold. 
	Then for Problem (P) and any $N$,
	\begin{equation}\label{eq13b}
	\max_{1\leq i\leq N} \mathbb{E}\int_0^{\infty} e^{-\rho t} \left(\|\hat{x}_i(t)\|^2+\|\hat{u}_i(t)\|^2\right)dt<\infty.
	\end{equation}
\end{theorem}
\emph{Proof.} {
From {\textbf{A5)}}, there exists $c_0>0$ such that $\mathbb{E}\int_0^{\infty}e^{-\rho t}\|\bar{x}(t)\|^2dt\leq c_0$. By 
Schwarz inequality,
\begin{align}\label{19}
  &\mathbb{E}\int_0^{\infty}e^{-\rho t}\|s_i(t)\|^2dt=\mathbb{E}\int_0^{\infty}e^{-\rho t}\left\|\int_t^{\infty}e^{-({\bar{A}_i}-{\rho}I)( t-\tau)}\bar{x}(\tau)d\tau\right\|^2dt\cr
=&  \mathbb{E}\int_0^{\infty}e^{-\rho t}\left\|\int_t^{\infty}e^{-(\frac{\bar{A}_i}{2}-\frac{\rho}{4}I+\frac{\bar{A}_i}{2}-\frac{3\rho}{4} I)( t-\tau)}\bar{x}(\tau)d\tau\right\|^2dt\cr
  \leq& \mathbb{E}\int_0^{\infty}e^{-\rho t}\int_t^{\infty}\big\|e^{-(\frac{\bar{A}_i}{2}-\frac{\rho}{4}I)( t-\tau)}\big\|^2d\tau
  \cdot\int_t^{\infty}\big\|e^{-(\frac{\bar{A}_i}{2}-\frac{3\rho}{4}I)( t-\tau)}\big\|^2\|\bar{x}(\tau)\|^2d\tau dt\cr
  \leq &\mathbb{E}\int_0^{\infty}\int_t^{\infty}\big\|e^{-(\frac{\bar{A}_i}{2}-\frac{\rho}{4}I)( t-\tau)}\big\|^2d\tau
  \cdot\int_t^{\infty}\big\|e^{-(\frac{\bar{A}_i}{2}-\frac{\rho}{4}I)( t-\tau)}\big\|^2e^{-\rho \tau}\|\bar{x}(\tau)\|^2d\tau dt.
\end{align}
By Lemma {\ref{thm1}}, $[\bar{A}_i,\bar{C}_i] $ is $\rho$-stable, which implies that all the eigenvalues of $\bar{A}_i-\frac{\rho}{2}I$
have negative real parts. Thus, there exist $c>0$ and $\delta_0>0$ such that  $\|e^{(\frac{\bar{A}_i}{2}-\frac{\rho}{4}I)t}\|^2\leq ce^{-\delta_0 t}.$
Note that $ \mathbb{E}\int_t^{\infty}e^{-\delta_0 (\tau-t)}d\tau=1/\delta_0. $ 
By (\ref{19}) and the exchange of order of the integration, we obtain
\begin{align*}
  \mathbb{E}\int_0^{\infty}e^{-\rho t}\|s_i(t)\|^2dt
  \leq &\mathbb{E}\int_0^{\infty}\frac{c}{\delta_0}
\int_t^{\infty}\big\|e^{-(\frac{\bar{A}_i}{2}-\frac{\rho}{4}I)( t-\tau)}\big\|^2e^{-\rho \tau}\|\bar{x}(\tau)\|^2d\tau dt\cr
=&\frac{c}{\delta_0}\mathbb{E}\int_0^{\infty}e^{-\rho \tau}\|\bar{x}(\tau)\|^2
\int_0^{\tau}\big\|e^{-(\frac{\bar{A}_i}{2}-\frac{\rho}{4}I)( t-\tau)}\big\|^2dtd\tau\cr
\leq &\frac{c^2}{\delta_0^2}\mathbb{E}\int_0^{\infty}(1-e^{-\delta_0 \tau})e^{-\rho \tau}\|\bar{x}(\tau)\|^2d\tau<\infty.
\end{align*}
By Lemma \ref{lem1} and \textbf{A1)},
\begin{equation*}\label{eq17b}
\begin{aligned}
\mathbb{E}\int_0^{\infty}e^{-\rho t}\|\hat{x}_i(t)\|^2dt\leq c\mathbb{E}\|x_i(0)\|^2+c\mathbb{E}\int_0^{\infty}e^{-\rho t}\|s_i(t)\|^2dt< \infty.
\end{aligned}
\end{equation*}
This together with (\ref{eq10bb}) completes the proof.}
\hfill$\Box$

\begin{remark}
	Theorem \ref{thm2} shows that the closed-loop system is uniformly $\rho$-stable. Particularly, when $\rho=0$, one can obtain that the closed-loop states of all the agents are asymptotically stable in the mean square sense.
\end{remark}

Define
\begin{equation}\label{eq25a}
  \epsilon_N^2=\int_0^{\infty} e^{-\rho t}\left\|\int_{\Theta}\bar{x}_{\theta}(t)dF_{N}(\theta)-\int_{\Theta}\bar{x}_{\theta}(t)dF(\theta)\right\|^2dt.
  \end{equation}

\begin{lemma}\label{thm3}
	Let \textbf{A1)}-\textbf{A5)} hold. Then we have
	$$
	\mathbb{ E}\int_0^{\infty} e^{-\rho t} \left\|\hat{x}^{(N)}(t)-\bar{x}(t)\right\|^2dt\leq c(\frac{1}{N}+\varepsilon_N^2),$$
	where $\hat{x}^{(N)} (t)=\frac{1}{N}\sum_{i=1}^N\hat{x}_i(t)$ and
	$\lim_{N\to \infty} \epsilon_N=0$.

\end{lemma}
\emph{Proof.} See Appendix \ref{app1}. \hfill{$\Box$}

\Blue{We now show the set of decentralized strategies is an $\varepsilon$-Nash equilibrium, whose definition is copied from \cite[Chapter 4]{BO82}.
\begin{definition}
For a given $\varepsilon\geq 0$, a set of strategies $\{u_i,\ 1 \leq i\leq N\}$ is called an $\varepsilon$-Nash equilibrium with respect to the set of cost functions
$\{J_i, 1 \leq i \leq N\}$ if for any $i=1, \cdots,  N$,
$$J_i(u_i, u_{-i}) \leq \inf_{u_i\in \mathcal{U}_{c}}
J_i(u_i, u_{-i}) + \varepsilon.$$
\end{definition}}

\begin{theorem}\label{thm4}
	Let \textbf{A1)}-\textbf{A5)} hold. Then for Problem (P), the set of strategies $(\hat{u}_1,\cdots,\hat{u}_N)$ given by (\ref{eq10bb}) is an $\varepsilon$-Nash equilibrium, i.e.,
	\begin{equation}
	J_i(\hat{u}_i, \hat{u}_{-i})-\varepsilon\leq \inf_{u_i\in {\mathcal U}_c}J_i({u}_i, \hat{u}_{-i})\leq J_i(\hat{u}_i, \hat{u}_{-i}),
	\end{equation}
	where $\varepsilon=O(\epsilon_N+\frac{1}{\sqrt{N}})$ and $\epsilon_N$ is given by (\ref{eq25a}).
\end{theorem}

To prove Theorem \ref{thm4}, we need the following lemma, whose proof is given in Appendix
B.

\begin{lemma}\label{lem2}
	Let \textbf{A1)}-\textbf{A5)} hold. Then
	\begin{align}
	J_i(\hat{u}_i, \hat{u}_{-i})\leq &\bar{J}_i(\hat{u}_i)+ \varepsilon,\label{eq17-a}
	\\
	\bar{J}_i(\hat{u}_i) \leq  & \inf_{u_i\in {\mathcal U}_c}J_i({u}_i, \hat{u}_{-i})+\varepsilon, \label{eq17}
	\end{align}
	where $\varepsilon=O(\epsilon_N+\frac{1}{\sqrt{N}})$.
\end{lemma}

\emph{Proof of Theorem \ref{thm4}.}  It follows from Lemma \ref{lem2} that
$$ J_i(\hat{u}_i, \hat{u}_{-i})\leq \inf_{u_i\in {\mathcal U}_c}J_i({u}_i, \hat{u}_{-i})+2\varepsilon,$$
which leads to the conclusion.   \hfill$\Box$

\section{The integrator systems with multiplicative noise}\label{Sec:multi-agent}

In this section, we first consider the single-integrator multi-agent systems with multiplicative noise, and then extend to the multiple-integrator systems.

\subsection{The model of noisy single integrator}
Now consider the degenerate case of (\ref{eq1}), where $x_i, u_i,W_i\in \mathbb{R}$, $A(\theta_i)=C=0$ and $B=D=Q=1$. In this situation, the dynamics of agent $i$ reduces to
\begin{equation}\label{eq12}
dx_i(t) = u_i(t)dt+u_i(t)dW_i(t), \ i=1,\cdots,N,
\end{equation}
The cost function of agent $i$ is given by
\begin{align}\label{eq11}
{J}_i(u_i,u_{-i})
=
\mathbb{E}\int_0^{\infty}e^{-\rho t}\big\{\big|x_i(t)
- {x}^{(N)}(t)\big|^2
+r|u_i(t)|^2\big\}dt,
\end{align}
\Blue{where $\rho>0$ and $r$ is allowed to negative.}

\textbf{Problem (P2).} \emph{For (\ref{eq12})-(\ref{eq11}), find an $\varepsilon$-Nash  equilibrium in ${\mathcal U}_{1}^{\prime}\times {\mathcal U}_{2}^{\prime}\times \cdots,{\mathcal U}_{N}^{\prime}$ to minimize $\{J_i,i=1,\cdots,N\}$, where
	\begin{equation*}
	\mathcal{U}_{i}^{\prime}=\Big\{u_i|u_i(t)\in \sigma(x_i(s),0\leq s\leq t), \mathbb{E}\int_0^{\infty} e^{-\rho t}|u_i(t)|^2dt<\infty
	\Big\}.
	\end{equation*}
}
We first give a preliminary result about the dynamics (\ref{eq12}).
\begin{lemma}\label{lem3}
	For (\ref{eq12}), there exists a constant $c_0$ such that
	$$\mathbb{E}\int_0^{\infty} e^{-\rho t}x_i^2(t)dt\leq
	c_0 \mathbb{E} \int_0^{\infty} e^{-\rho t}u_i^2(t)dt+c_0 , \ i=1,\cdots, N.$$
\end{lemma}
\emph{Proof.} Denote $x_{i,\rho}(t)=e^{-\frac{\rho}{2}t}x_i(t)$, and $u_{i,\rho}(t)=e^{-\frac{\rho}{2}t}u_i(t)$. It follows from (\ref{eq12}) that
$$dx_{i,\rho}(t)=-\frac{\rho}{2}x_{i,\rho}(t)dt+u_{i,\rho}(t)dt+u_{i,\rho}(t)dW_i(t).$$
By Lemma \ref{lem1}, we obtain the conclusion of this lemma. \hfill{$\Box$}

Replacing $x^{(N)}$ in (\ref{eq12}) by $\bar{x}\in C_{\rho/2}([0,\infty), \mathbb{R})$, we have
the limiting optimal control problem.

\textbf{Problem (P3).}  \emph{ Find $\bar{u}_i\in {\mathcal U}_{i}^{\prime}$ such that
	\begin{equation*}
	\bar{J}_i(\bar{u}_i)=\inf_{{u}_i\in {\mathcal U}_{i}^{\prime}}(u_i),
	\end{equation*} where
	$$\bar{J}_i(u_i)=\mathbb{E}\int_0^{\infty}e^{-\rho t}\Big\{\big|x_i(t)
	- \bar{x}(t)\big|^2
	+r|u_i(t)|^2\Big\}dt.$$
}
\begin{lemma}\label{lem3a}
	Let \textbf{A1)} hold and $r> \frac{2\sqrt{\rho+1}-(\rho+2)}{\rho^2}$. Suppose $\bar{x}\in C_{\rho/2}([0,\infty), \mathbb{R}^n)$ is given. For Problem (P3), the following statements hold.
	
	i) The algebraic Riccati equation
	\begin{equation}\label{eq11ab}
	\rho p=1-\frac{p^2}{r+p}
	\end{equation}
	admits a solution \begin{equation}\label{eq11ba}
	p=\frac{1-\rho r+\sqrt{\Delta}}{2(\rho+1)}
	\end{equation}
	such that $p+r>0$ and the system
	\begin{equation}\label{eq12a}
	dx_i(t)=-\frac{p}{p+r} x_i(t)dt-\frac{p}{p+r} x_i(t)dW_i(t)
	\end{equation}
	is $\rho$-stable. Here  $\Delta:=(\rho r)^2+2\rho r+4r+1$.

	ii) The optimal control law is uniquely given by
	$$\bar{u}_i(t)=-\frac{1}{p+r}(p x_i(t)+s(t)),\ t\geq0,$$
	where $s\in C_{\rho/2}([0,\infty), \mathbb{R}^n)$ is a unique solution to the equation
	\begin{equation}\label{eq12b}
	\frac{ds}{dt}= \big(\rho+\frac{p}{p+r}\big) s(t)+\bar{x}(t).
	\end{equation}

	iii) The optimal cost is given by
$$
		\inf_{u_i\in {\mathcal U}_{i}} \bar{J}_i(u_i)=
		p\mathbb{E}\|x_i(0)\|^2+2s(0)\bar{x}_0+g(0),
	$$
	where $g\in C_{\rho}([0,\infty), \mathbb{R})$ is a unique solution to the equation
		\begin{equation}\label{35a}
\frac{dg}{dt}=\rho g(t)+\frac{s^2(t)}{p+r}-\bar{x}^2(t).
\end{equation}
\end{lemma}

\emph{Proof.}  See Appendix B.    \hfill $\Box$

\begin{remark}\Blue{
In Lemma \ref{lem3a}, the condition $r> \frac{2\sqrt{\rho+1}-(\rho+2)}{\rho^2}$ is to ensure the convexity of Problem (P3). Note that the right hand is
always negative for $\rho > 0$, which implies that this condition is satisfied necessarily if $r > 0$.}
\end{remark}

For Problem (P2), the consistency equations (\ref{eq9a})-(\ref{eq11a}) reduce to
\begin{align}\label{eq14}
\frac{ds}{dt}=& \big(\rho+\frac{p}{p+r}\big) s(t)+\bar{x}(t),\\\label{eq15}
\frac{d\bar{x}}{dt}=&-\frac{1}{p+r}s(t) -\frac{p}{p+r}\bar{x}(t).
\end{align}

\begin{proposition}\label{thm5}
	Equations (\ref{eq14})-(\ref{eq15}) have a unique bounded solution
	\begin{equation*}
	(s(t),\bar{x}(t))\equiv (-p\bar{x}_0,\bar{x}_0), \quad t\geq0.
	\end{equation*}
\end{proposition}

\emph{Proof.} The equations (\ref{eq14})-(\ref{eq15}) can be written as
\begin{equation}\label{eq16}
\frac{d}{dt}\left[
\begin{array}{c}
s\\
\bar{x}
\end{array}\right]=H\left[
\begin{array}{c}
s(t)\\
\bar{x}(t)
\end{array}\right],
\end{equation}
where
$$H=\left[
\begin{array}{cc}
\rho+\frac{p}{p+r}& 1\\
-\frac{1}{p+r}&-\frac{p}{p+r}
\end{array}\right].$$
By some straightforward computations, we can diagonalize the matrix $H$ via $J=T^{-1}HT$, where
$$J= \left[
\begin{array}{cc}
0& 0\\
0&\rho
\end{array}\right], \quad T= \left[
\begin{array}{cc}
-p& p+r\\
1&-p
\end{array}\right].$$
The solution to (\ref{eq14})-(\ref{eq15}) given by
$$
\left[
\begin{array}{c}
s(t)\\
\bar{x}(t)
\end{array}\right]=e^{Ht}\left[
\begin{array}{c}
s(0)\\
\bar{x}(0)
\end{array}\right]$$  is bounded if and only if
$(s(0),\bar{x}(0))^T$ lies in the characteristic subspace of the eigenvalue 0. This implies $(s(0),\bar{x}(0))^T=c_0(-p,1)^T,$
where $c_0$ is a constant. Thus, we obtain that $s(0)=-p \bar{x}_0$ and (\ref{eq14})-(\ref{eq15}) has a unique solution
$$\left[
\begin{array}{c}
s(t)\\
\bar{x}(t)
\end{array}\right]=e^{Ht}\left[
\begin{array}{c}
-p \bar{x}_0\\
\bar{x}_0
\end{array}\right] \equiv\left[
\begin{array}{c}
-p \bar{x}_0\\
\bar{x}_0
\end{array}\right].$$
\rightline{$\Box$}

For a finite-population system in Problem (P2), we may construct the decentralized strategies as follows:
\begin{equation}\label{eq19}
\hat{u}_i(t)=-\frac{p}{p+r}(x_i(t)-\bar{x}_0), \quad t\geq 0, \quad i=1,\cdots, N.
\end{equation}
Substituting (\ref{eq19}) into (\ref{eq12}), the closed-loop dynamics of agent $i$ can be written as
$$d\hat{x}_i(t)=-\frac{p}{p+r}(\hat{x}_i(t)-\bar{x}_0)dt-\frac{p}{p+r}(\hat{x}_i(t)-\bar{x}_0)dW_i(t).$$
By It\^{o}'s formula, we have
\begin{equation}\label{eq20}
\begin{aligned}
\hat{x}_i(t)=\bar{x}_0+({x}_i(0)-\bar{x}_0)\exp &\Big[-\frac{p}{p+r}t-\frac{p^2}{2(p+r)^2}t-\frac{p}{p+r}W_i(t)\Big].
\end{aligned}
\end{equation}

\begin{lemma}
	\label{thm7}
	For Problem (P2), assume that \textbf{A1)} holds and  $r> \frac{2\sqrt{\rho+1}-(\rho+2)}{\rho^2}$. Under the set of strategies $\{\hat{u}_i,i=1,\cdots,N\}$ in (\ref{eq19}), we have
	\begin{equation}\label{eq21b}
	\lim_{N\to\infty}\int_0^{\infty}e^{-\rho t}\mathbb{E}|\hat{x}^{(N)}(t)-\bar{x}_0|^2dt=0.
	\end{equation}
\end{lemma}
\emph{Proof.} It follows from (\ref{eq20}) that
\begin{align*}
	\mathbb{ E}\Big|\frac{1}{N}\sum_{i=1}^N\hat{x}_i(t)-\bar{x}_0\Big|^2
	& =\ \mathbb{E}\left| \frac{1}{N}\sum_{i=1}^N({x}_i(0)-\bar{x}_0)\exp\Big[-\frac{p}{p+r}t -\frac{p^2}{2(p+r)^2}t-\frac{p}{p+r}W_i(t)\Big] \right|^2\\
	\hphantom{\mathbb{E}\Big|\frac{1}{N}\sum_{i=1}^N\hat{x}_i(t)-\bar{x}_0\Big|^2}
	&=\ \frac{1}{N^2}\sum_{i=1}^N\mathbb{E}|{x}_i(0)-\bar{x}_0|^2\mathbb{E}\exp\Big[-\frac{2p}{p+r}t-\frac{p^2}{(p+r)^2}t-\frac{2p}{p+r}W_i(t)\Big]\\
	\hphantom{\mathbb{ E}\Big|\frac{1}{N}\sum_{i=1}^N\hat{x}_i(t)-\bar{x}_0\Big|^2}
	&\leq \ \frac{1}{N}\sup_{1\leq i\leq N} \mathbb{E}|{x}_i(0)|^2\exp\Big[\Big(-\frac{2p}{p+r}+\frac{p^2}{(p+r)^2}\Big)t\Big].\\
\end{align*}
This gives
\begin{equation}\label{eq34}
\begin{aligned}
\int_0^{\infty}e^{-\rho t}\mathbb{E}|\hat{x}^{(N)}(t)-\bar{x}_0|^2dt =\frac{1}{N}\sup_{1\leq i\leq N} \mathbb{E}|{x}_i(0)|^2 \int_0^{\infty}\exp\Big[\Big(-\rho -\frac{2p}{p+r}+\frac{p^2}{(p+r)^2}\Big)t\Big]dt.
\end{aligned}
\end{equation}
By (\ref{eq11ab}) and (\ref{eq11ba}), we have
\begin{equation}\label{eq35}
-\rho -\frac{2p}{p+r}+\frac{p^2}{(p+r)^2}=\frac{1-\rho r-2(\rho+1)p}{p+r}=-\frac{\sqrt{\Delta}}{p+r},
\end{equation}
where $$\Delta={(\rho r)^2+2\rho r+4r+1}.$$ Note that
$$\frac{2\sqrt{\rho+1}-(\rho+2)}{\rho^2}=-\frac{1}{2\sqrt{\rho+1}+(\rho+2)}.$$ It can be verified that $\Delta>0$ when $r> \frac{2\sqrt{\rho+1}-(\rho+2)}{\rho^2}$.
By (\ref{eq11ba}),
$$p+r=\frac{(\rho+2)r+1+\sqrt{\Delta}}{2(\rho+1)}> \frac{1-\frac{(\rho+2)}{2\sqrt{\rho+1}+\rho+2}}{2(\rho+1)}>0.$$
From (\ref{eq34}) and  (\ref{eq35}), we obtain (\ref{eq21b}).   \hfill $\Box$

Now we give the result of asymptotic optimality. Denote $${\mathcal U}_{c}^{\prime}=\left\{u_i|u_i(t)\in \sigma(\bigcup_{i=1}^N {\mathcal F}_t^i),\mathbb{E}\int_0^{\infty} e^{-\rho t}|u_i(t)|^2dt<\infty\right\}.$$

\begin{theorem}\label{thm8}
	For Problem (P2), let \textbf{A1)} hold and  $r> \frac{2\sqrt{\rho+1}-(\rho+2)}{\rho^2}$. Then the set of strategies $(\hat{u}_1,\cdots,\hat{u}_N)$ given by (\ref{eq19}) is an $\varepsilon$-Nash equilibrium, i.e.,
	\begin{equation}\label{eq42}
	J_i(\hat{u}_i, \hat{u}_{-i})-\varepsilon\leq \inf_{u_i\in {\mathcal U}_c^{\prime}}J_i({u}_i, \hat{u}_{-i})\leq J_i(\hat{u}_i, \hat{u}_{-i}),
	\end{equation}
	where $\varepsilon=O(\epsilon_N+\frac{1}{\sqrt{N}})$. Furthermore, the asymptotic optimal cost of agent $i$ is $p(\mathbb{E}x_i^2(0)-\bar{x}_0^2)$.
\end{theorem}

To show Theorem \ref{thm8}, we first provide two lemmas, whose proofs are given in Appendix \ref{app3}.

\begin{lemma}\label{lem4}
Assume that \textbf{A1)} holds and  $r> \frac{2\sqrt{\rho+1}-(\rho+2)}{\rho^2}$.
Then 
there exists $\delta>0$ and $c>0$ such that
	\begin{align*}
	\bar{J}_i(u_i)
	&\geq \delta \mathbb{E}\int_0^{\infty}e^{-\rho t}u_i(t)dt-c.
	\end{align*}
\end{lemma}


\begin{lemma}\label{lem5}
	For Problem (P2), assume that \textbf{A1)} holds and  $r> \frac{2\sqrt{\rho+1}-(\rho+2)}{\rho^2}$. For $u_i\in {\mathcal U}_c$, if $J(u_i, \hat{u}_{-i})\leq c_1$, then there exist an integer $N_0$ and a constant $c_2$  such that for all $N\geq N_0$, $\mathbb{E}\int_0^{\infty} e^{-\rho t}u_i^2(t)dt \leq  c_2$.
\end{lemma}


\begin{remark}
Lemma \ref{lem4} shows under \textbf{A1)} and $r> \frac{2\sqrt{\rho+1}-(\rho+2)}{\rho^2}$, Problem (P3) is coercive. From this together with Lemma \ref{lem3a}, Problem (P3) is convex and coercive. The properties ensure the existence and asymptotic optimality of the decentralized strategies.
\end{remark}

\emph{Proof of Theorem \ref{thm8}.} Note that $\inf_{u_i\in {\mathcal U}_c}J_i({u}_i, \hat{u}_{-i})\leq J_i(\hat{u}_i, \hat{u}_{-i}).$ It suffices to show the second inequality in (\ref{eq42}) by checking all
$u_i$ such that
\begin{equation}
J_i({u}_i, \hat{u}_{-i})\leq J_i(\hat{u}_i, \hat{u}_{-i})\leq c_0.
\end{equation}
By Lemma \ref{lem5}, there exists an integer $N_0$ and a constant $c_2$  such that for all $N\geq N_0$, $$\mathbb{E}\int_0^{\infty} e^{-\rho t}u_i^2(t)dt\leq  c_2,$$
which together with Lemma \ref{lem3} gives
$$\mathbb{E}\int_0^{\infty} e^{-\rho t}x_i^2(t)dt\leq  c.$$
Following an argument similar to the proof of Lemma \ref{lem2}, we can obtain the second inequality in (\ref{eq42}).
By Lemma \ref{lem3a} and (\ref{eq11ab}), the asymptotic optimal cost is given by
$$\begin{aligned}
\lim_{N\to\infty}J_i(\hat{u}_i, \hat{u}_{-i})=\bar{J}_i(\hat{u}_i)
= \mathbb{E}[px_i^2(0)+2s(0)x_i(0)+g(0)] =p(\mathbb{E}x_i^2(0)-\bar{x}_0^2).
\end{aligned}$$
This completes the proof. \hfill $\Box$

We now show under mild conditions all the agents achieve 
mean-square consensus.

\begin{definition}\label{def1}	
	In a multi-agent system, all the agents  are said to reach
	the mean-square consensus if there exists a random variable
	$x^*$ such that
	$ \lim_{t\to \infty}\mathbb{E}|x_i(t)-x^*|^2=0$. 
	\end{definition}

%

\begin{theorem}\label{thm6}
	For Problem (P2), assume that \textbf{A1)} holds and $r>-\frac{1}{2(2+\rho)}$. Then there exist $c_1,c_2>0$ such that under the set of strategies $\{\hat{u}_i,i=1,\cdots,N\}$ in (\ref{eq19}),
	\begin{equation}\label{eq21}
	\mathbb{ E}|\hat{x}_i(t)-\bar{x}_0|^2\leq c_1e^{-c_2t}, \quad t>0, \quad 1\leq i\leq N,
	\end{equation}
	which imply
	all the agents achieve mean-square consensus.
\end{theorem}
\emph{Proof.}  It follows from (\ref{eq20}) that
\begin{align}\label{eq22}
\mathbb{E}|\hat{x}_i(t)-\bar{x}_0|^2 &=\mathbb{E}|{x}_i(0)-\bar{x}_0|^2\mathbb{E}\exp\Big[-\frac{2p}{p+r}t-\frac{p^2}{(p+r)^2}t-\frac{2p}{p+r}W_i(t)\Big] \nonumber \\
&= \mathbb{E}|{x}_i(0)-\bar{x}_0|^2\exp\Big[\Big(-\frac{2p}{p+r}+\frac{p^2}{(p+r)^2}\Big)t\Big] \nonumber \\
&\leq \sup_{1\leq i\leq N} \mathbb{E}|{x}_i(0)|^2\exp\Big[-\frac{p(p+2r)}{(p+r)^2}t\Big].
\end{align}
If $r\geq 0$, then we directly have $p>0$ and further obtain (\ref{eq21}). For the case $-\frac{1}{2(2+\rho)}< r< 0$, we first check the sign of $p+2r$.
By (\ref{eq11ba}),
\begin{align*}
p+2r=&\frac{(3\rho+4)r+1+\sqrt{(\rho r+1)^2+4r}}{2(\rho+1)}\cr
> & \frac{-\frac{\rho}{2(\rho+2)}+\sqrt{(1-\frac{\rho}{2(\rho+2)})^2-\frac{2}{\rho+2}}}{2(\rho+1)}\cr
=&\frac{-\rho+\rho}{4(\rho+1)(\rho+2)}=0.
\end{align*}
In view of $-\frac{1}{2(2+\rho)}< r< 0$, it follows that $p>p+r>p+2r>0$.
This together with (\ref{eq22}) implies (\ref{eq21}).
By \textbf{A1)} and the law of large numbers, 
 we have $\lim_{N\to \infty}\mathbb{E}|\frac{1}{N}\sum_{i=1}^N\hat{x}_i(0)-\bar{x}_0|^2=0$.
\hfill $\Box$

\begin{remark}\Blue{
Note that $2(\rho+2)-(\rho+2+2\sqrt{\rho+1})=(\sqrt{\rho+1}-1)^2>0$.
Then we have $$\frac{2\sqrt{\rho+1}-(\rho+2)}{\rho^2}=-\frac{1}{2\sqrt{\rho+1}+(\rho+2)}<-\frac{1}{2(2+\rho)}.$$
It means that $r>-\frac{1}{2(2+\rho)}$ implies $r> \frac{2\sqrt{\rho+1}-(\rho+2)}{\rho^2}$. Thus, to achieve mean-square consensus for agents, a tighter condition for $r$
is needed.}
\end{remark}

\subsection{The model of noisy multiple integrator}
We now extend the above result to the high-dimensional case. Specifically, agent $i$ evolves by
\begin{equation}\label{eq419}
  dx_i=Bu_idt+Du_idW_i,\ i=1,\cdots,N,
\end{equation}
and the cost function is given by
\begin{align}\label{eq420}
{J}_i(u_i,u_{-i})
=
\mathbb{E}\int_0^{\infty}e^{-\rho t}\big\{\big\|x_i(t)
- {x}^{(N)}(t)\big\|^2
+\|u_i(t)\|^2_R\big\}dt,
\end{align}
{where $\rho\geq0$ and $R$ is symmetric.}

For this problem, consistency equations (\ref{eq19a})-(\ref{eq19b}) have the following solution.
\begin{proposition}
Equations (\ref{eq19a})-(\ref{eq19b}) admit a unique bounded solution
	$
	(s(t),\bar{x}(t))\equiv (-P\bar{x}_0,\bar{x}_0), \ t\geq0.
	$
\end{proposition}
\emph{Proof.} Let $\Pi=P+K$,
where $P$ is the unique stabilizing solution to (\ref{eq6}) (i.e., a solution such that $A-BR^{-1}B^TP-\frac{\rho}{2}I$ is Hurwitz) and $K$ is the unique stabilizing solution to (\ref{eq24}). Note that $A=0$. We have that $\Pi$ satisfies
\begin{equation}\label{eq47}
  \rho \Pi+\Pi B(R+D^TPD)B^T\Pi=0.
\end{equation}
From (\ref{eq26}), ${A}-\frac{\rho}{2} I-B(R+D^TPD)^{-1}B^T\Pi$ is Hurwitz. Thus, the solution to (\ref{eq47}) is $\Pi=0$. In view of the proof of Theorem \ref{thm3.2}, we obtain $\bar{x}(t)\equiv \bar{x}_0$,
and $s(t)\equiv-P\bar{x}_0$. \hfill$\Box$

For a finite-population system in the above problem, we may construct the decentralized strategies as follows:
\begin{equation}\label{eq19m}
\hat{u}_i(t)=-(R+D^TPD)^{-1}B^T{P}(x_i(t)-\bar{x}_0), \quad t\geq 0, \quad i=1,\cdots, N.
\end{equation}
Substituting (\ref{eq19m}) into (\ref{eq12}), the closed-loop dynamics of agent $i$ can be written as
\begin{equation}\label{eq423}
  d\hat{x}_i(t)=-B(R+D^TPD)^{-1}B^T{P}(\hat{x}_i(t)-\bar{x}_0)dt-D(R+D^TPD)^{-1}B^T{P}(\hat{x}_i(t)-\bar{x}_0)dW_i(t).
  \end{equation}
\begin{theorem}
 The agents in (\ref{eq423}) reach mean-square consensus if and only if $[-B(R+D^TPD)^{-1}B^T{P},-D(R+D^TPD)^{-1}B^T{P}]$ is stable. Particularly, if $R>0$ and $\rho=0$,
  the agents achieve mean-square consensus.
\end{theorem}
\emph{Proof.}  Note that (\ref{eq423}) is equivalent to
$$d(\hat{x}_i(t)-\bar{x}_0)=-B(R+D^TPD)^{-1}B^T{P}(\hat{x}_i(t)-\bar{x}_0)dt-D(R+D^TPD)^{-1}B^T{P}(\hat{x}_i(t)-\bar{x}_0)dW_i(t).$$
From Definition \ref{def1}, we obtain that (\ref{eq423}) reaches mean-square consensus if and only if $[-B(R+D^TPD)^{-1}B^T{P},-D(R+D^TPD)^{-1}B^T{P}]$ is stable.
If $R>0$ and $\rho=0$, then by Lemma \ref{thm1}, $[-B(R+D^TPD)^{-1}B^T{P},-D(R+D^TPD)^{-1}B^T{P}]$ is stable. The theorem follows. \hfill$\Box$
\section{Numerical example}\label{Sec:example}
In this section, we give a numerical example to illustrate the effectiveness of the proposed consensus strategies for multi-agent systems with multiplicative noises. Furthermore, we compare the consensus effect between this case and multi-agent systems with additive noises \cite{NCMH13}, \cite{NCM14}.

First, consider a system of 50 agents with dynamics (\ref{eq1})-(\ref{eq2}), where $A(\theta_i) = 0.1$, $C = 0.1$, $B = D = 1$, $r = 1$ and $\rho = 0.6$. The initial states of the agents are taken independently from a normal distribution N(1, 0.5). After the strategies in (16) are applied, the state trajectories of 50 agents are shown in Figure 1. It can be seen that all the agents reach an agreement, but their states do not converge to a constant.

\begin{center}
	\vskip 0.3cm $\epsfig{figure=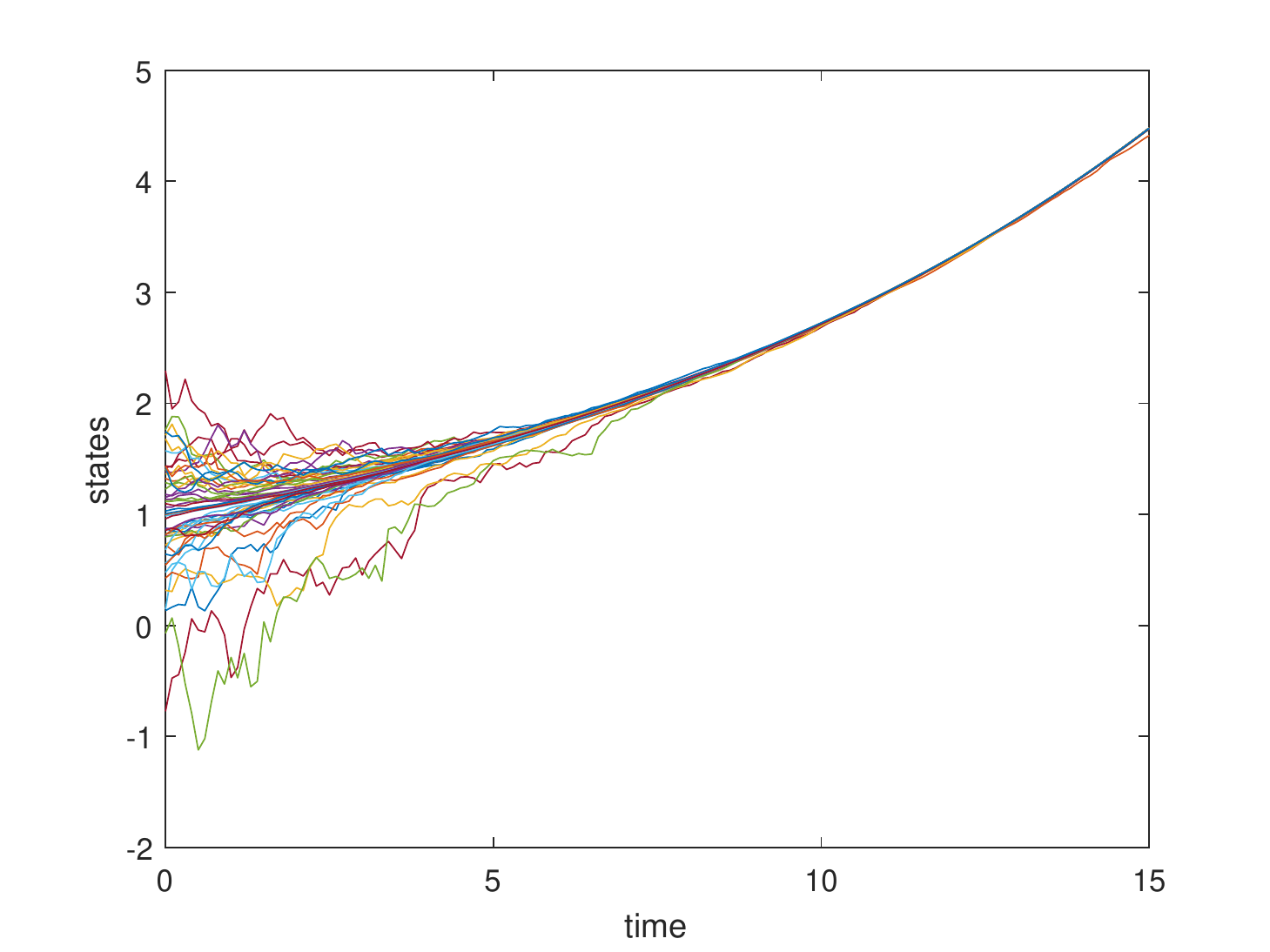,height=5.5cm}$\\
	{\small Figure 1: Under multiplicative noises 50 agents reach consensus ($A=0.1$)}
\end{center}

Next, consider a system of 50 agents with the single-integrator dynamics in (\ref{eq12})-(\ref{eq11}), where $D=0.5$, $r=1$ and $\rho=0.6$. The initial states of the agents are taken the same as above. Under strategies (\ref{eq19}), the state trajectories of the agents are shown in Figure 2. It is evident that the agents reach mean-square consensus.

\begin{center}
	\vskip 0.3cm $\epsfig{figure=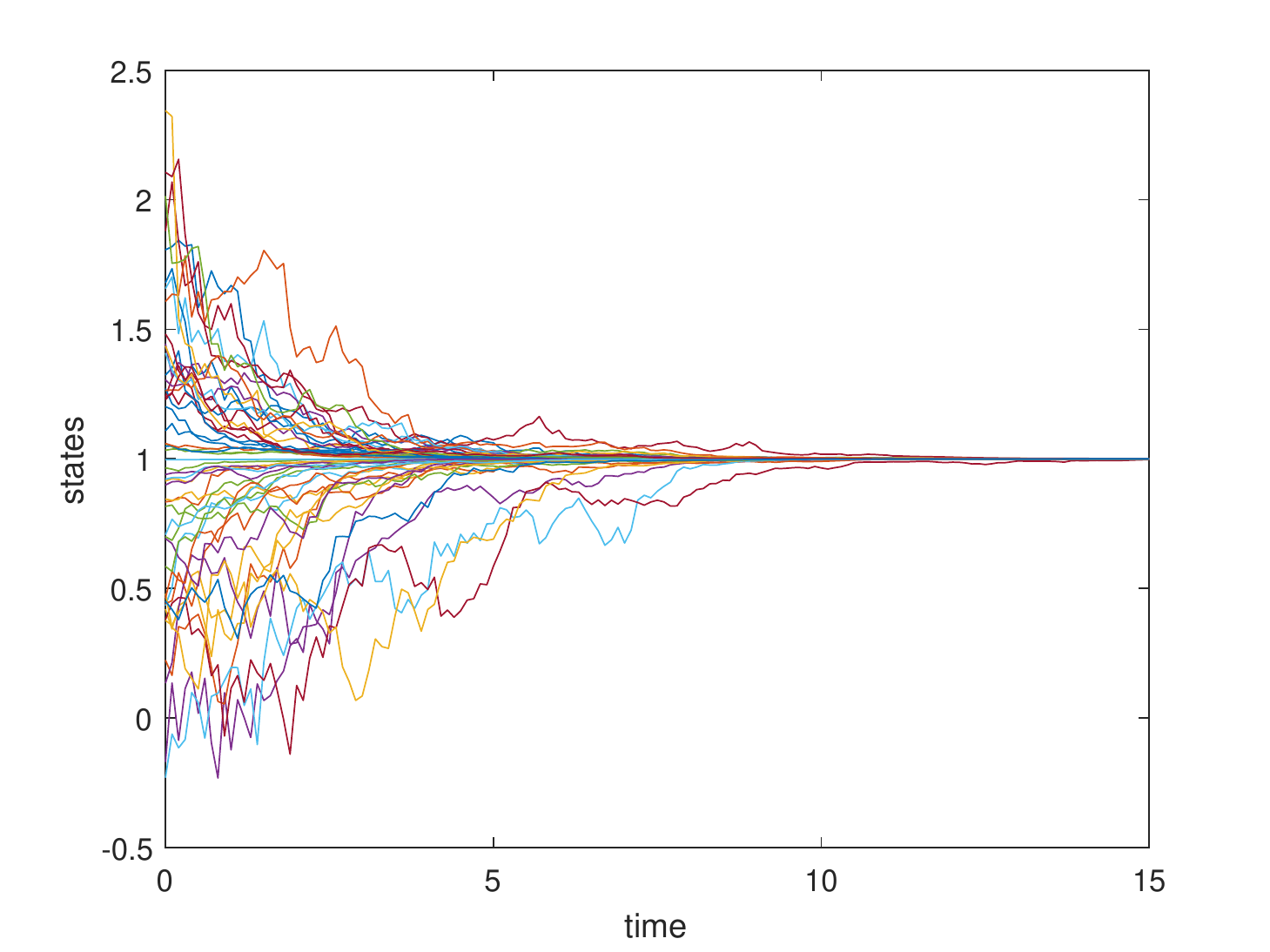,height=5.5cm}$\\
	{\small Figure 2: Under multiplicative noises 50 single-integrator agents reach consensus}
\end{center}

For comparison, we simulate the single-integrator multi-agent system with additive noise. Specifically, the dynamics of agent $i$ is given by
$$dx_i(t)=u_i(t)dt+\sigma dW_i(t),$$ and the cost function is $${J}_i(u_i,u_{-i})=\mathbb{E}\int_0^{\infty}e^{-\rho t}\big\{\big|x_i(t)
- {x}^{(N)}(t)\big|^2
+r|u_i(t)|^2\big\}dt.$$
The parameters $r, \rho$ are taken the same as above and the noise intensity $\sigma$ is $0.5$. The control strategy is designed by the mean field game methodology as in \cite{NCMH13}, \cite{NCM14}.
The state trajectories of the single-integrator multi-agent systems with additive noise are shown in Figure 3. It can be seen that the agents do not reach consensus, although they behave similarly.

\begin{center}
	\vskip 0.3cm $\epsfig{figure=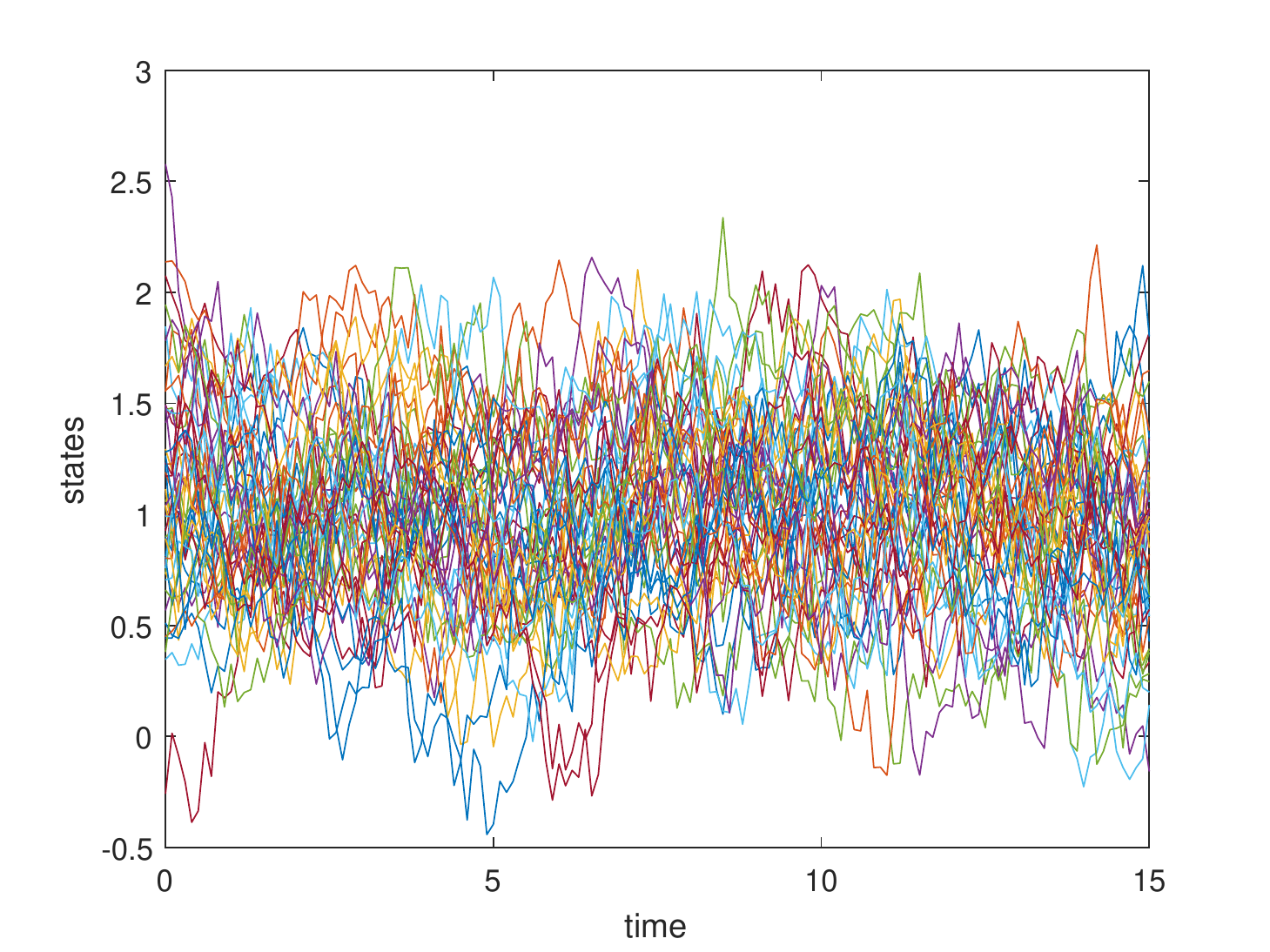,height=6cm}$\\
	{\small Figure 3: Under additive  noises 50 single-integrator agents do not reach consensus}
\end{center}


\section{Concluding remarks}\label{Sec:conclude}

In this paper, mean field games have been studied for multi-agent systems with multiplicative noises. By solving an auxiliary limiting optimal control problems subject to
consistent mean field approximations, a set of decentralized strategies is obtained, which is shown to be an $\varepsilon$-Nash equilibrium. The result is further applied to the consensus problem for the integrator multi-agent systems with multiplicative noises.
It is shown that under mild conditions the agents achieve asymptotic mean-square consensus.
For further works, it is of interest to consider mean field games for multi-agent systems with model uncertainty and common noises.

\bibliographystyle{plain}

\appendix

\section{Proofs of Lemma \ref{thm1} and Lemma \ref{thm3}}\label{app1}
\def\theequation{A.\arabic{equation}}
\setcounter{equation}{0}

\emph{Proof of Lemma \ref{thm1}.}
 i) Let $x_{i,\rho}(t)=e^{-\frac{\rho}{2}t}x_i(t)$, and $u_{i,\rho}(t)=e^{-\frac{\rho}{2}t}u_i(t)$. From (\ref{eq1}), we have
\begin{equation}\label{10b}
\begin{aligned}
dx_{i,\rho}(t)=	[(A_i-\frac{\rho}{2}I)x_{i,\rho}(t)+Bu_{i,\rho}](t)dt+[Cx_{i,\rho}(t)+Du_{i,\rho}(t)]dW_i(t).
\end{aligned}
\end{equation}
By \textbf{A4)} and \cite{ZZC08}, 
the equation (\ref{eq6}) admits a unique semi-positive definite solution, and under the control
$$\bar{u}_{i,\rho}(t)=-(R+D^TP_iD)^{-1}(B^TP_i+D^TP_iC) x_{i,\rho}(t),\ t\geq 0,$$
the solution of closed-loop system satisfies
$\mathbb{E}\int_0^{\infty}\|x_{i, \rho}(t)\|^2dt<\infty$. 
 Then, part i) follows.

ii) Letting $y_{i,\rho}(t)=\mathbb{E}[x_{i,\rho}(t)]$, we have from (\ref{10b})
\begin{equation*}
dy_{i,\rho}(t)= (A_i-\frac{\rho}{2}I)y_{i,\rho}(t)dt+Bu_{i,\rho}(t)dt.
\end{equation*}
Under the control
$\bar{u}_{i,\rho}=-(R+D^TP_iD)^{-1}(B^TP_i+D^TP_iC) x_{i,\rho}$, it holds that
$$ \int_0^{\infty}\|y_{i,\rho}(t)\|^2dt=\int_0^{\infty}\|\mathbb{E}x_{i,\rho}(t)\|^2dt\leq \mathbb{E}\int_0^{\infty}\|x_{i, \rho}(t)\|^2dt<\infty.$$
This implies that the real part of any eigenvalue of $\bar{A}_i-(\rho/2)I$ is negative. It follows from (\ref{eq7a}) that
$$s_i(t)=e^{-(\bar{A}_i-\rho I)^Tt}\Big[s_i(0)+\int_0^te^{(\bar{A}_i-\rho I)^T\tau}Q\bar{x}(\tau)d\tau\Big].$$
From an argument in \cite[Appendix A]{WZ12}, we have that $s\in C_{\rho/2}([0,\infty),\mathbb{R}^n)$ if and only if
the initial state $s(0)$ has the following expression
$
s_i(0)=-\int_0^{\infty}e^{(\bar{A}_i-\rho I)^T\tau}Q\bar{x}(\tau)d\tau.
$
Under this initial condition, we have
\begin{equation}\label{11a}
s_i(t)=-\int_t^{\infty}e^{-(\bar{A}_i-\rho I)^T(t-\tau)}Q\bar{x}(\tau)d\tau.
\end{equation}
Furthermore, using a similar argument as the above, we have $g_i\in C_{\rho}([0,\infty), \mathbb{R}^2)$ if and only if
\begin{align*}
  g_i(0)
  =-\int_0^{\infty}e^{-\rho t}
[s_i^T(t)B(R+D^TP_iD)^{-1}B^Ts_i(t)-\bar{x}^T(t)Q\bar{x}(t)]dt.
\end{align*}

iii) Applying It\^{o}'s formula to $e^{-\rho t}[x_i^T(t)P_ix_i(t)+2s_i^T(t)x_i(t)+g_i(t)]$, we have
\begin{eqnarray}\label{eq43}
&&\bar{J}_i(u_i)=\mathbb{E}[x_i^T(0)P_ix_i(0)+2s_i^T(0)x_i(0)+g_i(0)]-\lim_{T\to \infty}\mathbb{E}\big[ e^{-\rho T}(x_i^T(T)P_ix_i(T)+2s_i^T(T)x_i(T)+g_i(T))\big]\nonumber \\
&&\hphantom{\bar{J}_i(u_i)=}  +\mathbb{E}\int_0^{\infty}e^{-\rho t}\big\|u_i(t)+(R+D^TP_iD)^{-1}[(B^TP_i+D^TP_iC) x_i(t)+B^Ts_i(t)]\big\|^2_{R+D^TP_iD}dt\nonumber \\
&&\hphantom{\bar{J}_i(u_i)}=\mathbb{E}[x_i^T(0)P_ix_i(0)+2s_i^T(0)x_i(0)+g_i(0)] \cr
&&\hphantom{\bar{J}_i(u_i)=}+\mathbb{E}\int_0^{\infty}e^{-\rho t}\big\|u_i(t)+(R+D^TP_iD)^{-1}[(B^TP_i+D^TP_iC) x_i(t)+B^Ts_i(t)]\big\|^2_{R+D^TP_iD}dt\nonumber \\
&&\hphantom{\bar{J}_i(u_i)}\geq \mathbb{E}[x_i^T(0)P_ix_i(0)+2s_i^T(0)x_i(0)+g_i(0)].
\end{eqnarray}
Thus, (\ref{optimal-L}) is the optimal control and the optimal value is
$$\begin{aligned}
\bar{J}_i(\bar{u}_i)&=\inf_{u_i\in {\mathcal U}_{i}} \bar{J}_i(u_i)=\mathbb{E}[x_i^T(0)P_ix_i(0)+2s_i^T(0)x_i(0)+g_i(0)].
\end{aligned}
$$
This completes the proof.  \hfill$\Box$

\emph{Proof of Lemma \ref{thm3}.}
Note that $\{\hat{x}_i(t)-\bar{x}_i(t), t\geq 0\}, i=1,...,N$, are mutually independent processes with zero mean, where $\bar{x}_i(t)\stackrel{\Delta}{=}
\mathbb{E}x_i(t)=\bar{x}_{\theta_i}(t)$. We have
\begin{align*}
\mathbb{ E}\|\hat{x}^{(N)}-\bar{x}\|^2=&\mathbb{E}\Big\|\frac{1}{N}\sum_{i=1}^N(\hat{x}_i-\bar{x}_i)+\frac{1}{N}\sum_{i=1}^N(\bar{x}_i-\bar{x})\Big\|^2\cr
=&\frac{1}{N^2}\sum_{i=1}^N\mathbb{E}\|\hat{x}_i-\bar{x}_i\|^2+\Big\|\frac{1}{N}\sum_{i=1}^N\bar{x}_i-\bar{x}\Big\|^2.
\end{align*}
Note that $\hat{x}_i-\bar{x}_i$ satisfies
$$d(\hat{x}_i-\bar{x}_i)=\bar{A}_i(\hat{x}_i-\bar{x}_i)dt+\bar{C}\hat{x}_idW_i.$$
It follows from Theorem \ref{thm2} that
\begin{equation}\label{eq13a}
\begin{aligned}
&\frac{1}{N^2}\mathbb{E}\int_0^{\infty}e^{-\rho t}\sum_{i=1}^N\|\hat{x}_i-\bar{x}_i\|^2dt \cr
\leq&\frac{1}{N}\max_{1\leq i\leq N}\mathbb{E}\int_0^{\infty}\big\|e^{(\bar{A}_i-
\frac{\rho}{2} I)t}({x}_{i0}-\bar{x}_0)\big\|^2dt+ \frac{1}{N}\max_{1\leq i\leq N}\mathbb{E}\int_0^{\infty}e^{-\rho t}\Big\|\int_0^te^{\bar{A}_i(t-\tau )}\hat{x}_i(\tau )dW_i(\tau)\Big\|^2dt\cr
\leq& \frac{1}{N}\int_0^{\infty}\big\|e^{(\bar{A}_i-
\frac{\rho}{2} I)t}\big\|^2\max_{1\leq i\leq N}\mathbb{E}\|{x}_{i0}-\bar{x}_0)\|^2dt+\frac{1}{N}\max_{1\leq i\leq N}\mathbb{E}\int_0^{\infty}e^{-\rho t}\int_0^t\big\|e^{\bar{A}_i(t-\tau )}\big\|^2\|\hat{x}_i(\tau )\|^2d\tau dt
\cr
\leq &\frac{c_1}{N}+\frac{1}{N}\max_{1\leq i\leq N}\mathbb{E}\int_0^{\infty}e^{-\rho \tau}\|\hat{x}_i(\tau )\|^2\int_{\tau}^{\infty}\big\|e^{(\bar{A}_i-\frac{\rho}{2}I)(t-\tau )}\big\|^2 dtd\tau
\leq \frac{c_2}{N}.
\end{aligned}
\end{equation}
By \textbf{A3)}, we have
$$\left\|\frac{1}{N}\sum_{i=1}^N\bar{x}_i-\bar{x}\right\|=\left\|\int_{\Theta}\bar{x}_{\theta}dF_{N}(\theta)-\int_{\Theta}\bar{x}_{\theta}dF(\theta)\right\|\stackrel{\Delta}{=} h_N(t).$$
By the weak convergence in \textbf{A2)}, $\lim_{N\to\infty}h_N(t)=0$ for each $t$, which further implies $$\lim_{N\to\infty}\epsilon^2_N=\lim_{N\to\infty}\int_0^{\infty} e^{-\rho t}h_N(t)dt=0.$$
Recalling (\ref{eq13a}), the lemma follows. \hfill$\Box$

\section{Proofs of Lemma \ref{lem2} and Lemma \ref{lem3a}}
\def\theequation{B.\arabic{equation}}
\setcounter{equation}{0}
\emph{Proof of Lemma \ref{lem2}.} By Cauchy-Schwarz inequality and Lemma \ref{thm3}, we have
\begin{eqnarray*}
&&J_i(\hat{u}_i, \hat{u}_{-i})=\mathbb{E}\int_0^{\infty}e^{-\rho t}\Big\{\big\|\hat{x}_i(t)-\bar{x}(t)+\bar{x}(t)-\hat{x}^{(N)}(t)\big\|^2_{Q} +\|\hat{u}_i(t)\|^2_{R}\Big\}dt\\
&&\hphantom{J_i(\hat{u}_i, \hat{u}_{-i})}\leq \bar{J}_i(\hat{u}_i)+\mathbb{E}\int_0^{\infty}e^{-\rho t}\big\|\bar{x}-\hat{x}^{(N)}\big\|^2_{Q}dt + 2c\left(\mathbb{E}\int_0^{\infty}e^{-\rho t}\big\|\bar{x}-\hat{x}^{(N)}\big\|^2_{Q}dt\right)^{1/2}\\
&&\hphantom{J_i(\hat{u}_i, \hat{u}_{-i})}\leq  \bar{J}_i(\hat{u}_i)+ \varepsilon.
\end{eqnarray*}
Note that $\inf_{u_i\in {\mathcal U}_c}J_i({u}_i, \hat{u}_{-i})\leq J_i(\hat{u}_i, \hat{u}_{-i}).$ It suffices to show (\ref{eq17}) by checking all
$u_i$ such that
\begin{equation}\label{eq18}
J_i({u}_i, \hat{u}_{-i})\leq J_i(\hat{u}_i, \hat{u}_{-i})\leq c_0.
\end{equation}
It follows by Theorem \ref{thm2} that
$$\begin{aligned}
\mathbb{E}\int_0^{\infty}e^{-\rho t}\big\|{x}_i
- \frac{1}{N}{x}_i\big\|^2_{Q}dt&\leq \mathbb{E}\int_0^{\infty}e^{-\rho t}\big\|{x}_i
- \frac{1}{N}{x}_i-\frac{1}{N}\sum_{j\not=i}\hat{x}_j\big\|^2_{Q}dt +\mathbb{E}\int_0^{\infty}e^{-\rho t}\big\|\frac{1}{N}\sum_{j\not=i}\hat{x}_j\big\|^2_{Q}dt
\cr&\leq c_0+ \mathbb{E}\int_0^{\infty}e^{-\rho t}\frac{1}{N-1}\sum_{j\not=i}\big\|\hat{x}_j\big\|^2_{Q}dt\leq c,
\end{aligned}$$
which implies
\begin{equation}\label{eq18-b}
\mathbb{E}\int_0^{\infty}e^{-\rho t}\|{x}_i
\|^2_{Q}dt\leq c_1,
\end{equation}
where $c_1$ is independent of $N$. By direct calculation, we have
\begin{align}\label{eq18-c}
J_i({u}_i, \hat{u}_{-i}) &=\mathbb{E}\int_0^{\infty}e^{-\rho t}\Big\{\big\|{x}_i
- \bar{x}+\bar{x}-\frac{1}{N}x_i -\frac{1}{N}\sum_{j\not=i}\hat{x}_j\big\|^2_{Q}
+\|\hat{u}_i\|^2_{R}\Big\}dt\cr
&\geq  \bar{J}_i(u_i)+2\mathbb{E}\int_0^{\infty}e^{-\rho t}({x}_i
- \bar{x})^TQ\Big(\bar{x}-\frac{1}{N}x_i-\frac{1}{N}\sum_{j\not=i}\hat{x}_j\Big)dt\cr
&\stackrel{\Delta}{=}\bar{J}_i(u_i)+2I_1.
\end{align}
By Cauchy-Schwarz inequality and (\ref{eq18-b}),
\begin{align*}
I^2_1\leq& \mathbb{E}\int_0^{\infty}e^{-\rho t}\big(2\|{x}_i
\|^2_Q+\| \bar{x}\|^2_Q\big)dt  \mathbb{E}\int_0^{\infty}e^{-\rho t}  \Big(2\big\|\frac{1}{N}\bar{x}-\frac{1}{N}x_i\big\|^2_Q+2\big\|\bar{x}-\hat{x}^{(N)}\big\|^2_Q\Big)\cr
\leq &c(\frac{c}{N^2}+2\varepsilon^2).
\end{align*}
From this together with (\ref{eq18-c}), we get (\ref{eq17}).  \hfill$\Box$

\emph{Proof of Lemma \ref{lem3a}.}  i)
By solving (\ref{eq11ab}),  we obtain
\begin{equation*}
p=\frac{1-\rho r\pm\sqrt{\Delta}}{2(\rho+1)},
\end{equation*}
where  $\Delta=(\rho r)^2+2\rho r+4r+1$.
Choose the maximal solution \cite{RZ00}
\begin{equation}\label{eq11b}
p=\frac{1-\rho r+\sqrt{\Delta}}{2(\rho+1)}.
\end{equation}
Note that $$\frac{2\sqrt{\rho+1}-(\rho+2)}{\rho^2}=-\frac{1}{2\sqrt{\rho+1}+(\rho+2)}<-\frac{1}{\rho}.$$ It can be verified that $\Delta>0$ when $r> \frac{2\sqrt{\rho+1}-(\rho+2)}{\rho^2}$.
By (\ref{eq11b}),
\begin{align*}
p+r
&>\frac{\frac{(\rho+2)[2\sqrt{\rho+1}-(\rho+2)]}{\rho^2}+1}{2(\rho+1)}=\frac{1-\frac{(\rho+2)}{2\sqrt{\rho+1}+\rho+2}}{2(\rho+1)}>0.
\end{align*}
It follows by (\ref{eq12a}) that
\begin{equation*}
\frac{dx_i}{x_i}=-\frac{p}{p+r} dt-\frac{p}{p+r} dW_i(t).
\end{equation*}
Applying It\^{o}'s formula to $\ln x_i$ and using the above equation, we have
\begin{align*}
d(\ln x_i)
= & \frac{dx_i}{x_i}-\frac{p^2}{2(p+r)^2}dt\cr
=& -\Big[\frac{p}{p+r}+\frac{p^2}{2(p+r)^2}\Big]dt-\frac{p}{p+r}dW_i(t),
\end{align*}
which implies
$$x_i(t)=x_i(0)\exp\Big[-\frac{p}{p+r}t-\frac{p^2}{2(p+r)^2}t-\frac{p}{p+r}W_i(t)\Big].$$
This follows a geometric Brownian motion.
By A1) and the direct calculation,
$$\mathbb{E}|{x}_i(t)|^2= \mathbb{E}|{x}_i(0)|^2\exp\Big[\Big(-\frac{2p}{p+r}+\frac{p^2}{(p+r)^2}\Big)t\Big],$$
which gives
\begin{equation*}
\begin{aligned}
&\int_0^{\infty}e^{-\rho t}\mathbb{E}|{x}_i(t)|^2dt=\mathbb{E}|{x}_i(0)|^2 \int_0^{\infty}\exp\Big[\Big(-\rho -\frac{2p}{p+r}+\frac{p^2}{(p+r)^2}\Big)t\Big]dt.
\end{aligned}
\end{equation*}
By (\ref{eq11ab}) and (\ref{eq11b}),
\begin{equation*}
\begin{aligned}
&-\rho -\frac{2p}{p+r}+\frac{p^2}{(p+r)^2}
=\frac{1-\rho r-2(\rho+1)p}{p+r}=-\frac{\sqrt{\Delta}}{p+r}<0.
\end{aligned}
\end{equation*}
Thus,  $ \int_0^{\infty}e^{-\rho t}\mathbb{E}|{x}_i(t)|^2dt<\infty$, i.e., (\ref{eq12a}) is $\rho$-stable.

ii)-iii) It follows from (\ref{eq12b}) that $$s(t)=e^{(\rho+\frac{p}{p+r})t}\Big[s(0)+\int_0^te^{-(\rho+\frac{p}{p+r})\tau}\bar{x}d\tau\Big].$$
Since $p+r>0$, it is straightforward to show that
$s\in C_{\rho/2}([0,\infty), \mathbb{R}^n)$ if and only if $s(0)=-\int_0^{\infty}e^{-(\rho+\frac{p}{p+r})\tau}\bar{x}d\tau,$
which implies
\begin{equation*}\label{eq12d}
s(t)=-\int_t^{\infty}e^{\big(\rho+\frac{p}{p+r}\big)(t-\tau)}\bar{x}d\tau.
\end{equation*}
By a similar argument, we have $g\in C_{\rho}([0,\infty), \mathbb{R}^n)$ if and only if
$$g(0)=-\int_0^{\infty}e^{-\rho\tau}\big[\frac{s^2}{p+r}-\bar{x}^2\big]d\tau.$$
By Lemma \ref{lem3}, for any $u_i\in {\mathcal{U}}_i^{\prime}$, $\mathbb{E}\int_0^{\infty}e^{-\rho t}x_i^2dt<\infty$.
Applying It\^{o}'s formula into $e^{-\rho t}(px_i^2+2sx_i+g)$, we have
\begin{align}\label{eq43}
\bar{J}_i(u_i)
=&\mathbb{E}[px_i^2(0)+2s(0)x_i(0)+g(0)] \cr
&+\mathbb{E}\int_0^{\infty}e^{-\rho t}(p+r)\Big[u_i+\frac{p}{p+r}(x_i-\bar{x})\Big]^2dt\cr
\geq& \mathbb{E}[px_i^2(0)+2s(0)x_i(0)+g(0)].
\end{align}
Thus, if and only if $u_i=\bar{u}_i=-\frac{p}{p+r}(x_i-\bar{x})$, then
$$\bar{J}_i(\bar{u}_i)=\inf_{u_i\in {\mathcal U}^{\prime}_{i}} \bar{J}_i(u_i)= \mathbb{E}[px_i^2(0)+2s(0)x_i(0)+g(0)].$$
This completes the proof. \hfill{$\Box$}

\section{Proofs of Lemma \ref{lem4} and Lemma \ref{lem5}}\label{app3}
\def\theequation{C.\arabic{equation}}
\setcounter{equation}{0}

\emph{Proof of Lemma \ref{lem4}. } From (\ref{35a}) and Proposition {\ref{thm5}},
it follows that $g$ satisfies
$$\frac{dg}{dt}=\rho g-\rho p\bar{x}_0^2.$$
We can identify this equation admits a unique solution $g(t)\equiv p\bar x_0^2$ in $C_{\rho}([0,\infty),\mathbb{R})$. By Lemma \ref{lem4}, $\mathbb{E}\int_0^{\infty} e^{-\rho t}x_i^2(t)dt\leq c.$
By It\^{o}'s formula, we have
\begin{align}\label{eq43}
\bar{J}_i(u_i)=&\mathbb{E}[px_i^2(0)+2s(0)x_i(0)+g(0)]-\lim_{T\to \infty}\mathbb{E}\big[ e^{-\rho T}(px_i^2(T)+2s(T)x_i(T)+g(T))\big]\cr
&+\mathbb{E}\int_0^{\infty}e^{-\rho t}(p+r)\Big[u_i+\frac{p}{p+r}(x_i-\bar{x}_0)\Big]^2dt\cr
=&\mathbb{E}[px_i^2(0)]-p\bar x_0^2  +\mathbb{E}\int_0^{\infty}e^{-\rho t}(p+r)\Big[u_i+\frac{p}{p+r}(x_i-\bar{x}_0)\Big]^2dt.
\end{align}
Let $u_i^*=u_i+\frac{p}{p+r}(x_i-\bar{x}_0)$. Then by (\ref{eq12}), $x_i$ satisfies
$$dx_i=\big[-\frac{p}{p+r}(x_i-\bar{x}_0)+u^*_i\big]dt+\big[-\frac{p}{p+r}(x_i-\bar{x}_0)+u^*_i\big]dW_i. $$
By Lemma \ref{lem1}, there exists a constant $c>0$,
$$\mathbb{E}\int_0^{\infty}e^{-\rho t}|x_i|^2dt\leq c_0 \mathbb{E}\int_0^{\infty}e^{-\rho t}|u_i^*|^2dt+c_0.$$
Thus,
\begin{align*}
  \mathbb{E}\int_0^{\infty}e^{-\rho t}|u_i|^2dt=&\mathbb{E}\int_0^{\infty}e^{-\rho t}[u_i^*-\frac{p}{p+r}(x_i-\bar{x}_0)]^2dt\cr
  \leq&c\mathbb{E}\int_0^{\infty}e^{-\rho t}|u_i^*|^2dt+c\cr
  =& c\mathbb{E}\int_0^{\infty}e^{-\rho t}\big[u_i+\frac{p}{p+r}(x_i-\bar{x}_0)\big]^2dt+c.
\end{align*}
Note that $p+r>0$.
From (\ref{eq43}), the lemma follows.
\hfill{$\Box$}

\emph{Proof of Lemma \ref{lem5}.} For $u_i\in {\mathcal U}_c$, we have
\begin{align}\label{eq45}
c_1\geq & J(u_i,\hat{u}_{-i})\cr
=&\mathbb{E}\int_0^{\infty}e^{-\rho t}\left\{\Big[x_i-\frac{1}{N}\Big(x_i+\sum_{j\not=i}^N\hat x_j\Big)\Big]^2+ru_i^2\right\}dt\cr
=&\mathbb{E}\int_0^{\infty}e^{-\rho t}\left\{\Big[x_i-\bar{x}+\bar{x}-\frac{1}{N}\Big(x_i+\sum_{j\not=i}^N\hat x_j\Big)\Big]^2+ru_i^2\right\}dt\cr
\geq &\bar{J}(u_i)+2 \mathbb{E}\int_0^{\infty}e^{-\rho t}(x_i-\bar{x})\Big[\bar{x}-\frac{1}{N}\Big(x_i+\sum_{j\not=i}^N\hat x_j\Big)\Big]dt\cr
\geq& \bar{J}(u_i)-2I_2,
\end{align}
where
\begin{equation*}
\begin{aligned}
I_2\stackrel{\Delta}{=}2&\left\{ \mathbb{E}\int_0^{\infty}e^{-\rho t}(x_i-\bar{x})^2dt \cdot \mathbb{E}\int_0^{\infty}e^{-\rho t}\Big[\bar{x}-\hat{x}^{(N)}-\frac{1}{N}(x_i-\hat x_i)\Big]^2dt\right\}^{\frac{1}{2}}.
\end{aligned}
\end{equation*}
Notice that $\bar{x}\equiv\bar{x}_0$. By Lemma \ref{lem3}, we have
\begin{equation}\label{eq46}
\begin{aligned}
\mathbb{E}\int_0^{\infty}e^{-\rho t}(x_i-\bar{x})^2dt &\leq  \mathbb{E}\int_0^{\infty}e^{-\rho t}(2x^2_i+2\bar{x}^2)dt \cr
&\leq c \mathbb{E}\int_0^{\infty} e^{-\rho t}u_i^2dt+c.
\end{aligned}
\end{equation}
It follows from Theorems \ref{thm7} and \ref{thm2} that
$$\begin{aligned}
&\mathbb{E}\int_0^{\infty}e^{-\rho t} \Big[\bar{x}-\hat{x}^{(N)}-\frac{1}{N}(x_i-\hat x_i)\Big]^2dt \cr
\leq & \mathbb{E}\int_0^{\infty}e^{-\rho t} \Big[2(\bar{x}-\hat{x}^{(N)})^2+\frac{2}{N^2}(x_i-\hat x_i)^2\Big]dt\cr
\leq& \frac{c}{N}+\frac{c}{N^2}  \mathbb{E}\int_0^{\infty} e^{-\rho t}u_i^2dt,
\end{aligned}
$$
which together with (\ref{eq45}) implies
$$|I_2|\leq \frac{c}{N}\mathbb{E}\int_0^{\infty} e^{-\rho t}u_i^2dt+\frac{c}{\sqrt{N}}. $$
By (\ref{eq45}) and Lemma \ref{lem4}, we obtain
$$c_1\geq J(u_i, \hat{u}_{-i})\geq(\delta-\frac{c}{N})\mathbb{E}\int_0^{\infty} e^{-\rho t}u_i^2dt-2c.$$
Letting $N_0=\inf\{m\in \mathbb{Z}|m>c/\delta\}$, we obtain that there exists a constant $c_2$ such that for all $N>N_0$,
$$\mathbb{E}\int_0^{\infty} e^{-\rho t}u_i^2(t)dt\leq\frac{N(c_1+2c)}{N\delta-c}\stackrel{\Delta}{=}c_2.$$ This completes the proof.
\hfill{$\Box$}


%

\end{document}